\documentclass[11pt]{article}
\usepackage{latexsym}

\usepackage[top=1in, bottom=1in, left=1in, right=1in]{geometry}

\usepackage{adjustbox}

\usepackage{amssymb}

\usepackage{graphicx} 
\usepackage{color} 
\usepackage{amsmath, amsthm, amssymb}
\usepackage{enumerate}
\usepackage{float}
\usepackage{url}
\usepackage{stackrel}
\usepackage{mathrsfs,dsfont}
\usepackage{caption}
\usepackage{subcaption}
\usepackage{wrapfig}
\usepackage{bbm}
\usepackage{bm}
\usepackage{epigraph}
\usepackage{physics}
\usepackage[colorinlistoftodos, shadow]{todonotes} 

\usepackage{mdwlist}

\usepackage{pifont}

\newtheorem{theorem}{Theorem}[section]

\newtheorem{example}[theorem]{Example}
\newtheorem{remark}[theorem]{Remark}

\theoremstyle{definition}
\newtheorem{definition}[theorem]{Definition}
\theoremstyle{definition}

\newcommand{\been}{\begin{enumerate}}
\newcommand{\enen}{\end{enumerate}}
\newcommand{\beit}{\begin{itemize}}
\newcommand{\enit}{\end{itemize}}

\def\SS{\mathcal S}
\def\CC{\mathcal C}

\def\PP{\mathcal P}

\def\GG{\mathcal{G}}

\def\HH{\mathcal{H}}

\def\DD{\mathcal{D}}

\def\spn{{\rm span}}

\def\ve{\varepsilon}

\def\la{\leftarrow}

\def\rlas{\rightleftarrows}
\def\ds{\displaystyle}

\def\cbl{\color{blue}}
\def\cma{\color{magenta}}
\def\ccy{\color{cyan}}
\def\cte{\color{teal}}

\newcommand{\R}{\mathbb{R}}

\DeclareMathOperator{\sgn}{sgn}

\usepackage{tikz}    
\usetikzlibrary{shapes,automata,positioning,arrows,fit}
\usepackage{mathtools}  
\usetikzlibrary{snakes}

\usepackage{caption}
\newcommand{\specialcell}[2][c]{\begin{tabular}[#1]{@{}c@{}}#2\end{tabular}}

\newcommand{\alis}[1]{\begin{align*}#1\end{align*}}

\numberwithin{equation}{section}

\usepackage{tikz-3dplot}

\usepackage{relsize}

 \tikzset{every node/.style={auto}}
 \tikzset{every state/.style={rectangle, minimum size=0pt, draw=none, font=\normalsize}}
  \tikzset{bend angle=7}
  
  \usepackage{array}
  
  \def\wt{\widetilde}

\usepackage{authblk}

\begin{document}

\title{Reaction Network Motifs for Static and Dynamic Absolute Concentration Robustness}

\author{Badal Joshi \footnote{Department of Mathematics, California State University San Marcos (bjoshi$@$csusm.edu)} ~ and Gheorghe Craciun \footnote{Departments of Mathematics and Biomolecular Chemistry, University of Wisconsin-Madison (craciun$@$wisc.edu)}}

\date{}



\maketitle

\begin{abstract}
Networks with {\em absolute concentration robustness} (ACR) have the property that a translation of a coordinate hyperplane either contains all steady states (static ACR) or attracts all trajectories (dynamic ACR). 
The implication for the underlying biological system is  robustness in the concentration of one of the species independent of the initial conditions as well as independent of the concentration of all other species.
Identifying network conditions for dynamic ACR is a challenging problem. 
We lay the groundwork in this paper by  studying small reaction networks, those with 2 reactions and 2 species. 
We give a complete classification by ACR properties of these minimal reaction networks.  
The dynamics is rich even within this simple setting.  
Insights obtained from this work will help illuminate the properties of more complex networks with dynamic ACR. 
~ 
\end{abstract}

%

%
%

\section{Introduction}

Absolute Concentration Robustness (ACR) was introduced in \cite{shinar2010structural} as the mathematical property that every positive steady state of a reaction system coincides in one of the coordinates. 
We call this property static ACR.
Assuming convergence to one of these positive steady states, static ACR ensures that at steady state, the value of one of the measured variables will be independent of the initial value. 
However, such convergence is not guaranteed even for the simplest systems. 
The long-term behavior of a system can lead to the system converging to the boundary, or to another attractor such as a limit cycle or not converging at all but diverging to infinity. 
We introduced {\em dynamic ACR} in \cite{joshi2021foundations} to account for the {\em global dynamics}, and not merely the location of steady states. 

A network condition for static ACR is found in \cite{shinar2010structural}: if the network deficiency is one and two non-terminal complexes differ in a single species, then the concentration of that species will have static ACR for all mass action reaction rate constants. 
Since analyzing global dynamics of a dynamical system is an enormously more complex task than determining the locations of steady states, obtaining network conditions for dynamic ACR is challenging.
We make some headway in this direction by studying small reaction networks and organizing our ideas carefully in this setting. 
In this paper, we study networks with two reactions and two species and make fine distinctions between global dynamical properties related to convergence to an ACR hyperplane.
This lays the groundwork for future study of network conditions for larger, more biochemically realistic networks. 

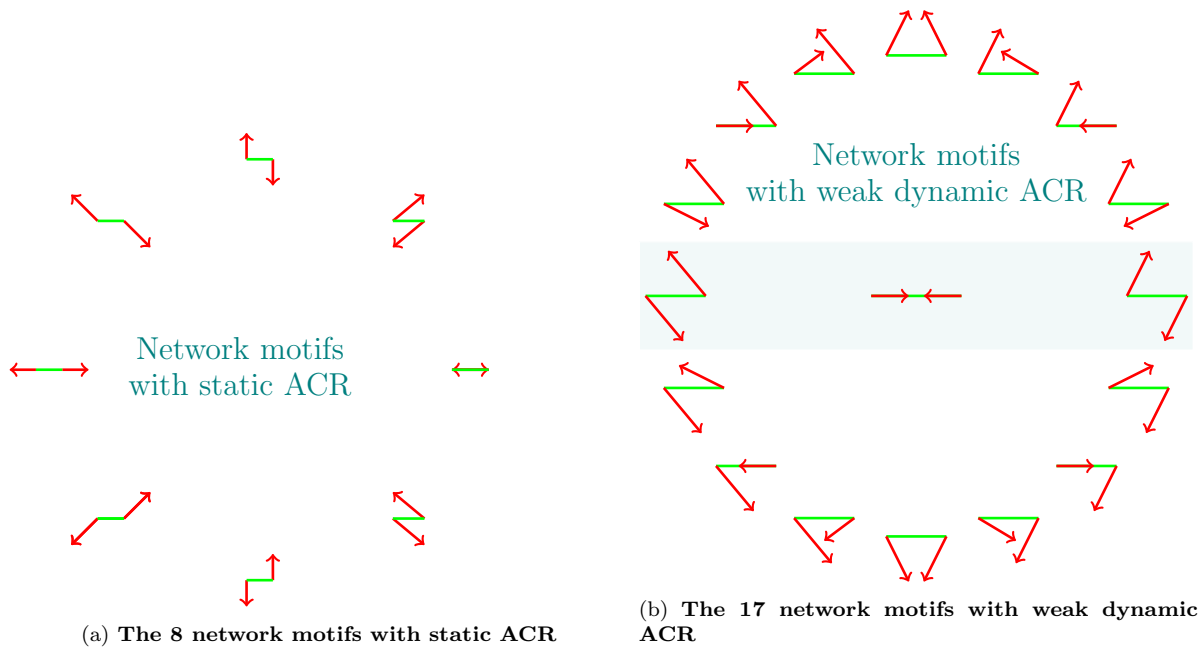
\begin{figure}[h!] 
\centering
\begin{subfigure}[b]{0.5\textwidth}
\begin{tikzpicture}[scale=0.35]
\draw [->, line width=1, red] ({8*cos(90)+0.5},{8*sin(90)}) -- ({8*cos(90)+0.5},{8*sin(90)-1});
\draw [->, line width=1, red] ({8*cos(90)-0.5},{8*sin(90)}) -- ({8*cos(90)-0.5},{8*sin(90)+1});
\draw [-, line width=1, green] ({8*cos(90)-0.5},{8*sin(90)}) -- ({8*cos(90)+0.5},{8*sin(90)});

\draw [->, line width=1, red] ({8*cos(45)+0.6},{8*sin(45)}) -- ({8*cos(45)-0.6},{8*sin(45)-1});
\draw [->, line width=1, red] ({8*cos(45)-0.6},{8*sin(45)}) -- ({8*cos(45)+0.6},{8*sin(45)+1});
\draw [-, line width=1, green] ({8*cos(45)-0.6},{8*sin(45)}) -- ({8*cos(45)+0.6},{8*sin(45)});

\draw [->, line width=1, red] ({8*cos(0)+0.7},{8*sin(0)}) -- ({8*cos(0)-0.7},{8*sin(0)});
\draw [->, line width=1, red] ({8*cos(0)-0.7},{8*sin(0)}) -- ({8*cos(0)+0.7},{8*sin(0)});
\draw [-, line width=1, green] ({8*cos(0)-0.7},{8*sin(0)}) -- ({8*cos(0)+0.7},{8*sin(0)});

\draw [->, line width=1, red] ({8*cos(-45)+0.6},{8*sin(-45)}) -- ({8*cos(-45)-0.6},{8*sin(-45)+1});
\draw [->, line width=1, red] ({8*cos(-45)-0.6},{8*sin(-45)}) -- ({8*cos(-45)+0.6},{8*sin(-45)-1});
\draw [-, line width=1, green] ({8*cos(-45)-0.6},{8*sin(-45)}) -- ({8*cos(-45)+0.6},{8*sin(-45)});

\draw [->, line width=1, red] ({8*cos(-90)+0.5},{8*sin(-90)}) -- ({8*cos(-90)+0.5},{8*sin(-90)+1});
\draw [->, line width=1, red] ({8*cos(-90)-0.5},{8*sin(-90)}) -- ({8*cos(-90)-0.5},{8*sin(-90)-1});
\draw [-, line width=1, green] ({8*cos(-90)-0.5},{8*sin(-90)}) -- ({8*cos(-90)+0.5},{8*sin(-90)});

\draw [->, line width=1, red] ({8*cos(-135)+0.5},{8*sin(-135)}) -- ({8*cos(-135)+1.5},{8*sin(-135)+1});
\draw [->, line width=1, red] ({8*cos(-135)-0.5},{8*sin(-135)}) -- ({8*cos(-135)-1.5},{8*sin(-135)-1});
\draw [-, line width=1, green] ({8*cos(-135)-0.5},{8*sin(-135)}) -- ({8*cos(-135)+0.5},{8*sin(-135)});

\draw [->, line width=1, red] ({8*cos(-135)+0.5},{8*sin(-135)}) -- ({8*cos(-135)+1.5},{8*sin(-135)+1});
\draw [->, line width=1, red] ({8*cos(-135)-0.5},{8*sin(-135)}) -- ({8*cos(-135)-1.5},{8*sin(-135)-1});
\draw [-, line width=1, green] ({8*cos(-135)-0.5},{8*sin(-135)}) -- ({8*cos(-135)+0.5},{8*sin(-135)});

\draw [->, line width=1, red] ({8*cos(180)+0.5},{8*sin(180)}) -- ({8*cos(180)+1.5},{8*sin(180)});
\draw [->, line width=1, red] ({8*cos(180)-0.5},{8*sin(180)}) -- ({8*cos(180)-1.5},{8*sin(180)});
\draw [-, line width=1, green] ({8*cos(180)-0.5},{8*sin(180)}) -- ({8*cos(180)+0.5},{8*sin(180)});

\draw [->, line width=1, red] ({8*cos(135)+0.5},{8*sin(135)}) -- ({8*cos(135)+1.5},{8*sin(135)-1});
\draw [->, line width=1, red] ({8*cos(135)-0.5},{8*sin(135)}) -- ({8*cos(135)-1.5},{8*sin(135)+1});
\draw [-, line width=1, green] ({8*cos(135)-0.5},{8*sin(135)}) -- ({8*cos(135)+0.5},{8*sin(135)});

\node[shape=ellipse, draw=none, line width=0]() at (-0.75,0) {\color{teal} \specialcell{\large Network motifs\\\large with static ACR}};
\end{tikzpicture}
\caption{{\bf The 8 network motifs with static ACR}}
\label{fig:qonqegl;nl;}
\end{subfigure}
\begin{subfigure}[b]{0.45\textwidth}
\begin{tikzpicture}[scale=0.4]

\filldraw[color=red!0, fill=teal!5] (-9.2,-1.8) rectangle (9.2,1.8);

\draw [-, line width=1, green] ({8*cos(90)-1},{8*sin(90)}) -- ({8*cos(90)+1},{8*sin(90)});
\draw [->, line width=1, red] ({8*cos(90)-1},{8*sin(90)}) -- ({8*cos(90)-0.25},{8*sin(90)+1.5});
\draw [->, line width=1, red] ({8*cos(90)+1},{8*sin(90)}) -- ({8*cos(90)+0.25},{8*sin(90)+1.5});

\draw [-, line width=1, green] ({8*cos(67.5)-1},{8*sin(67.5)}) -- ({8*cos(67.5)+1},{8*sin(67.5)});
\draw [->, line width=1, red] ({8*cos(67.5)-1},{8*sin(67.5)}) -- ({8*cos(67.5)-0.25},{8*sin(67.5)+1.5});
\draw [->, line width=1, red] ({8*cos(67.5)+1},{8*sin(67.5)}) -- ({8*cos(67.5)-0.25},{8*sin(67.5)+0.75});

\draw [-, line width=1, green] ({8*cos(45)-1},{8*sin(45)}) -- ({8*cos(45)+1},{8*sin(45)});
\draw [->, line width=1, red] ({8*cos(45)-1},{8*sin(45)}) -- ({8*cos(45)-0.25},{8*sin(45)+1.5});
\draw [->, line width=1, red] ({8*cos(45)+1},{8*sin(45)}) -- ({8*cos(45)-0.25},{8*sin(45)});

\draw [-, line width=1, green] ({8*cos(22.5)-1},{8*sin(22.5)}) -- ({8*cos(22.5)+1},{8*sin(22.5)});
\draw [->, line width=1, red] ({8*cos(22.5)-1},{8*sin(22.5)}) -- ({8*cos(22.5)-0.25},{8*sin(22.5)+1.5});
\draw [->, line width=1, red] ({8*cos(22.5)+1},{8*sin(22.5)}) -- ({8*cos(22.5)-0.5},{8*sin(22.5)-0.75});

\draw [-, line width=1, green] ({8*cos(0)-1},{8*sin(0)}) -- ({8*cos(0)+1},{8*sin(0)});
\draw [->, line width=1, red] ({8*cos(0)-1},{8*sin(0)}) -- ({8*cos(0)-0.25},{8*sin(0)+1.5});
\draw [->, line width=1, red] ({8*cos(0)+1},{8*sin(0)}) -- ({8*cos(0)+0.25},{8*sin(0)-1.5});

\draw [-, line width=1, green] ({8*cos(-22.5)-1},{8*sin(-22.5)}) -- ({8*cos(-22.5)+1},{8*sin(-22.5)});
\draw [->, line width=1, red] ({8*cos(-22.5)-1},{8*sin(-22.5)}) -- ({8*cos(-22.5)+0.5},{8*sin(-22.5)+0.75});
\draw [->, line width=1, red] ({8*cos(-22.5)+1},{8*sin(-22.5)}) -- ({8*cos(-22.5)+0.25},{8*sin(-22.5)-1.5});

\draw [-, line width=1, green] ({8*cos(-45)-1},{8*sin(-45)}) -- ({8*cos(-45)+1},{8*sin(-45)});
\draw [->, line width=1, red] ({8*cos(-45)-1},{8*sin(-45)}) -- ({8*cos(-45)+0.25},{8*sin(-45)});
\draw [->, line width=1, red] ({8*cos(-45)+1},{8*sin(-45)}) -- ({8*cos(-45)+0.25},{8*sin(-45)-1.5});

\draw [-, line width=1, green] ({8*cos(-67.5)-1},{8*sin(-67.5)}) -- ({8*cos(-67.5)+1},{8*sin(-67.5)});
\draw [->, line width=1, red] ({8*cos(-67.5)-1},{8*sin(-67.5)}) -- ({8*cos(-67.5)+0.25},{8*sin(-67.5)-0.75});
\draw [->, line width=1, red] ({8*cos(-67.5)+1},{8*sin(-67.5)}) -- ({8*cos(-67.5)+0.25},{8*sin(-67.5)-1.5});

\draw [-, line width=1, green] ({8*cos(-90)-1},{8*sin(-90)}) -- ({8*cos(-90)+1},{8*sin(-90)});
\draw [->, line width=1, red] ({8*cos(-90)-1},{8*sin(-90)}) -- ({8*cos(-90)-0.25},{8*sin(-90)-1.5});
\draw [->, line width=1, red] ({8*cos(-90)+1},{8*sin(-90)}) -- ({8*cos(-90)+0.25},{8*sin(-90)-1.5});

\draw [-, line width=1, green] ({8*cos(112.5)-1},{8*sin(112.5)}) -- ({8*cos(112.5)+1},{8*sin(112.5)});
\draw [->, line width=1, red] ({8*cos(112.5)-1},{8*sin(112.5)}) -- ({8*cos(112.5)-0},{8*sin(112.5)+0.75});
\draw [->, line width=1, red] ({8*cos(112.5)+1},{8*sin(112.5)}) -- ({8*cos(112.5)-0.25},{8*sin(112.5)+1.5});

\draw [-, line width=1, green] ({8*cos(135)-1},{8*sin(135)}) -- ({8*cos(135)+1},{8*sin(135)});
\draw [->, line width=1, red] ({8*cos(135)-1},{8*sin(135)}) -- ({8*cos(135)+0.25},{8*sin(135)});
\draw [->, line width=1, red] ({8*cos(135)+1},{8*sin(135)}) -- ({8*cos(135)-0.25},{8*sin(135)+1.5});

\draw [-, line width=1, green] ({8*cos(157.5)-1},{8*sin(157.5)}) -- ({8*cos(157.5)+1},{8*sin(157.5)});
\draw [->, line width=1, red] ({8*cos(157.5)-1},{8*sin(157.5)}) -- ({8*cos(157.5)+0.5},{8*sin(157.5)-0.75});
\draw [->, line width=1, red] ({8*cos(157.5)+1},{8*sin(157.5)}) -- ({8*cos(157.5)-0.25},{8*sin(157.5)+1.5});

\draw [-, line width=1, green] ({8*cos(180)-1},{8*sin(180)}) -- ({8*cos(180)+1},{8*sin(180)});
\draw [->, line width=1, red] ({8*cos(180)-1},{8*sin(180)}) -- ({8*cos(180)+0.25},{8*sin(180)-1.5});
\draw [->, line width=1, red] ({8*cos(180)+1},{8*sin(180)}) -- ({8*cos(180)-0.25},{8*sin(180)+1.5});

\draw [-, line width=1, green] ({8*cos(202.5)-1},{8*sin(202.5)}) -- ({8*cos(202.5)+1},{8*sin(202.5)});
\draw [->, line width=1, red] ({8*cos(202.5)-1},{8*sin(202.5)}) -- ({8*cos(202.5)+0.25},{8*sin(202.5)-1.5});
\draw [->, line width=1, red] ({8*cos(202.5)+1},{8*sin(202.5)}) -- ({8*cos(202.5)-0.5},{8*sin(202.5)+0.75});

\draw [-, line width=1, green] ({8*cos(225)-1},{8*sin(225)}) -- ({8*cos(225)+1},{8*sin(225)});
\draw [->, line width=1, red] ({8*cos(225)-1},{8*sin(225)}) -- ({8*cos(225)+0.25},{8*sin(225)-1.5});
\draw [->, line width=1, red] ({8*cos(225)+1},{8*sin(225)}) -- ({8*cos(225)-0.25},{8*sin(225)+0});

\draw [-, line width=1, green] ({8*cos(247.5)-1},{8*sin(247.5)}) -- ({8*cos(247.5)+1},{8*sin(247.5)});
\draw [->, line width=1, red] ({8*cos(247.5)-1},{8*sin(247.5)}) -- ({8*cos(247.5)+0.25},{8*sin(247.5)-1.5});
\draw [->, line width=1, red] ({8*cos(247.5)+1},{8*sin(247.5)}) -- ({8*cos(247.5)-0},{8*sin(247.5)-0.75});

\draw [-, line width=1, green] ({0-1.5},{0}) -- ({0+1.5},{0});
\draw [->, line width=1, red] ({0-1.5},{0}) -- ({0-0.25},{0});
\draw [->, line width=1, red] ({0+1.5},{0}) -- ({0+0.25},{0});

\node[shape=ellipse, draw=none, line width=0]() at ({4*cos(90)},{4*sin(90)}) {\color{teal} \specialcell{\large Network motifs \\\large with weak dynamic ACR}};
\end{tikzpicture}
\caption{{\bf The  17 network motifs with weak dynamic ACR}}
\end{subfigure}
\caption{{\bf Network motifs with output robustness:}  {\bf (Left)} All possible two-reaction two-species network motifs with static ACR.  {\bf (Right)} All possible two-reaction two-species network motifs for which the hyperplane $\{x=x^*\}$ is weakly stable (invariant and weakly attracting) (see Definition \ref{def:o4tyhoegjwppp}.) The motifs with one-dimensional stoichiometric space are on the horizontal band. The motifs on the circumference require at least two species while the motif in the center can be realized with one or two species.}
\label{fig:allthemotifs}
\end{figure}

The definition of dynamic ACR introduced in \cite{joshi2021foundations} is generalized here in two ways. The first is by including a basin of attraction for the ACR hyperplane, which may be the entire positive orthant or a proper subset of it. 
The second way is by considering a weaker version of attracting hyperplane: trajectories may not converge to the ACR hyperplane, but may only move in the direction of that hyperplane. 
For static ACR, we introduce a stronger form which requires at least one steady state in each compatibility class that intersects the ACR hyperplane. 
While each ACR property (static, strong static, dynamic, weak dynamic) is worthy of study by itself, there are connections between them which give a more complete picture. Weak dynamic ACR is a necessary condition for dynamic ACR; similarly, when a positive steady state exists, static ACR is necessary for dynamic ACR (see Theorem \ref{thm:formsofacr} and Figure \ref{fig:3oyho3thwef;w;o}). Therefore networks with static ACR and weak dynamic ACR serve as candidates for networks with dynamic ACR.
{\bf We show that all motifs of static ACR and weak dynamic ACR for networks with two reactions and two species can be completely characterized -- there are exactly 8 network motifs with static ACR and there are exactly 17 network motifs with weak dynamic ACR and these are depicted in Figure \ref{fig:allthemotifs}} (see Theorems \ref{thm:po43ghi3lwkfh}, \ref{thm:oh3ihgj3oth3} and \ref{thm:oerigheogh}; also see \cite{meshkat2021absolute} for a network characterization of static ACR).  

\begin{figure}[H] 
\centering
\begin{tikzpicture}[scale=0.6]

\filldraw[color=red!0, fill=teal!5] (-14,-1.8) rectangle (14,1.8);

\draw [-, line width=2, green] ({8*cos(90)-1},{8*sin(90)}) -- ({8*cos(90)+1},{8*sin(90)});
\draw [->, line width=2, red] ({8*cos(90)-1},{8*sin(90)}) -- ({8*cos(90)-0.25},{8*sin(90)+1.5});
\draw [->, line width=2, red] ({8*cos(90)+1},{8*sin(90)}) -- ({8*cos(90)+0.25},{8*sin(90)+1.5});

\draw [-, line width=2, green] ({8*cos(67.5)-1},{8*sin(67.5)}) -- ({8*cos(67.5)+1},{8*sin(67.5)});
\draw [->, line width=2, red] ({8*cos(67.5)-1},{8*sin(67.5)}) -- ({8*cos(67.5)-0.25},{8*sin(67.5)+1.5});
\draw [->, line width=2, red] ({8*cos(67.5)+1},{8*sin(67.5)}) -- ({8*cos(67.5)-0.25},{8*sin(67.5)+0.75});

\draw [-, line width=2, green] ({8*cos(45)-1},{8*sin(45)}) -- ({8*cos(45)+1},{8*sin(45)});
\draw [->, line width=2, red] ({8*cos(45)-1},{8*sin(45)}) -- ({8*cos(45)-0.25},{8*sin(45)+1.5});
\draw [->, line width=2, red] ({8*cos(45)+1},{8*sin(45)}) -- ({8*cos(45)-0.25},{8*sin(45)});

\draw [-, line width=2, green] ({8*cos(22.5)-1},{8*sin(22.5)}) -- ({8*cos(22.5)+1},{8*sin(22.5)});
\draw [->, line width=2, red] ({8*cos(22.5)-1},{8*sin(22.5)}) -- ({8*cos(22.5)-0.25},{8*sin(22.5)+1.5});
\draw [->, line width=2, red] ({8*cos(22.5)+1},{8*sin(22.5)}) -- ({8*cos(22.5)-0.5},{8*sin(22.5)-0.75});

\draw [-, line width=2, green] ({8*cos(0)-1},{8*sin(0)}) -- ({8*cos(0)+1},{8*sin(0)});
\draw [->, line width=2, red] ({8*cos(0)-1},{8*sin(0)}) -- ({8*cos(0)-0.25},{8*sin(0)+1.5});
\draw [->, line width=2, red] ({8*cos(0)+1},{8*sin(0)}) -- ({8*cos(0)+0.25},{8*sin(0)-1.5});

\draw [-, line width=2, green] ({8*cos(-22.5)-1},{8*sin(-22.5)}) -- ({8*cos(-22.5)+1},{8*sin(-22.5)});
\draw [->, line width=2, red] ({8*cos(-22.5)-1},{8*sin(-22.5)}) -- ({8*cos(-22.5)+0.5},{8*sin(-22.5)+0.75});
\draw [->, line width=2, red] ({8*cos(-22.5)+1},{8*sin(-22.5)}) -- ({8*cos(-22.5)+0.25},{8*sin(-22.5)-1.5});

\draw [-, line width=2, green] ({8*cos(-45)-1},{8*sin(-45)}) -- ({8*cos(-45)+1},{8*sin(-45)});
\draw [->, line width=2, red] ({8*cos(-45)-1},{8*sin(-45)}) -- ({8*cos(-45)+0.25},{8*sin(-45)});
\draw [->, line width=2, red] ({8*cos(-45)+1},{8*sin(-45)}) -- ({8*cos(-45)+0.25},{8*sin(-45)-1.5});

\draw [-, line width=2, green] ({8*cos(-67.5)-1},{8*sin(-67.5)}) -- ({8*cos(-67.5)+1},{8*sin(-67.5)});
\draw [->, line width=2, red] ({8*cos(-67.5)-1},{8*sin(-67.5)}) -- ({8*cos(-67.5)+0.25},{8*sin(-67.5)-0.75});
\draw [->, line width=2, red] ({8*cos(-67.5)+1},{8*sin(-67.5)}) -- ({8*cos(-67.5)+0.25},{8*sin(-67.5)-1.5});

\draw [-, line width=2, green] ({8*cos(-90)-1},{8*sin(-90)}) -- ({8*cos(-90)+1},{8*sin(-90)});
\draw [->, line width=2, red] ({8*cos(-90)-1},{8*sin(-90)}) -- ({8*cos(-90)-0.25},{8*sin(-90)-1.5});
\draw [->, line width=2, red] ({8*cos(-90)+1},{8*sin(-90)}) -- ({8*cos(-90)+0.25},{8*sin(-90)-1.5});

\draw [-, line width=2, green] ({8*cos(112.5)-1},{8*sin(112.5)}) -- ({8*cos(112.5)+1},{8*sin(112.5)});
\draw [->, line width=2, red] ({8*cos(112.5)-1},{8*sin(112.5)}) -- ({8*cos(112.5)-0},{8*sin(112.5)+0.75});
\draw [->, line width=2, red] ({8*cos(112.5)+1},{8*sin(112.5)}) -- ({8*cos(112.5)-0.25},{8*sin(112.5)+1.5});

\draw [-, line width=2, green] ({8*cos(135)-1},{8*sin(135)}) -- ({8*cos(135)+1},{8*sin(135)});
\draw [->, line width=2, red] ({8*cos(135)-1},{8*sin(135)}) -- ({8*cos(135)+0.25},{8*sin(135)});
\draw [->, line width=2, red] ({8*cos(135)+1},{8*sin(135)}) -- ({8*cos(135)-0.25},{8*sin(135)+1.5});

\draw [-, line width=2, green] ({8*cos(157.5)-1},{8*sin(157.5)}) -- ({8*cos(157.5)+1},{8*sin(157.5)});
\draw [->, line width=2, red] ({8*cos(157.5)-1},{8*sin(157.5)}) -- ({8*cos(157.5)+0.5},{8*sin(157.5)-0.75});
\draw [->, line width=2, red] ({8*cos(157.5)+1},{8*sin(157.5)}) -- ({8*cos(157.5)-0.25},{8*sin(157.5)+1.5});

\draw [-, line width=2, green] ({8*cos(180)-1},{8*sin(180)}) -- ({8*cos(180)+1},{8*sin(180)});
\draw [->, line width=2, red] ({8*cos(180)-1},{8*sin(180)}) -- ({8*cos(180)+0.25},{8*sin(180)-1.5});
\draw [->, line width=2, red] ({8*cos(180)+1},{8*sin(180)}) -- ({8*cos(180)-0.25},{8*sin(180)+1.5});

\draw [-, line width=2, green] ({8*cos(202.5)-1},{8*sin(202.5)}) -- ({8*cos(202.5)+1},{8*sin(202.5)});
\draw [->, line width=2, red] ({8*cos(202.5)-1},{8*sin(202.5)}) -- ({8*cos(202.5)+0.25},{8*sin(202.5)-1.5});
\draw [->, line width=2, red] ({8*cos(202.5)+1},{8*sin(202.5)}) -- ({8*cos(202.5)-0.5},{8*sin(202.5)+0.75});

\draw [-, line width=2, green] ({8*cos(225)-1},{8*sin(225)}) -- ({8*cos(225)+1},{8*sin(225)});
\draw [->, line width=2, red] ({8*cos(225)-1},{8*sin(225)}) -- ({8*cos(225)+0.25},{8*sin(225)-1.5});
\draw [->, line width=2, red] ({8*cos(225)+1},{8*sin(225)}) -- ({8*cos(225)-0.25},{8*sin(225)+0});

\draw [-, line width=2, green] ({8*cos(247.5)-1},{8*sin(247.5)}) -- ({8*cos(247.5)+1},{8*sin(247.5)});
\draw [->, line width=2, red] ({8*cos(247.5)-1},{8*sin(247.5)}) -- ({8*cos(247.5)+0.25},{8*sin(247.5)-1.5});
\draw [->, line width=2, red] ({8*cos(247.5)+1},{8*sin(247.5)}) -- ({8*cos(247.5)-0},{8*sin(247.5)-0.75});

\draw [-, line width=2, green] ({0-1.5},{0}) -- ({0+1.5},{0});
\draw [->, line width=2, red] ({0-1.5},{0}) -- ({0-0.25},{0});
\draw [->, line width=2, red] ({0+1.5},{0}) -- ({0+0.25},{0});

\node[shape=ellipse, draw=none, line width=0]() at ({1*cos(-90)},{1*sin(-90)}) {\color{cyan} \specialcell{Full basin\\DACR}};

\node[shape=ellipse, draw=none, line width=0]() at ({10.5*cos(90)},{10.5*sin(90)}) {\color{cyan} \specialcell{Full basin\\DACR}};

\node[shape=ellipse, draw=none, line width=0]() at ({11.5*cos(0)},{11.5*sin(0)}) {\specialcell{{\color{magenta}Neighborhood \&}\\{\color{magenta}Almost Cylinder}\\{\color{cyan}Narrow basin}\\{\color{cyan} full space DACR}}};

\node[shape=ellipse, draw=none, line width=0]() at ({11.5*cos(22.5)},{11.5*sin(22.5)})  {\specialcell{{\color{magenta}Cylinder DACR}\\{\color{cyan}Narrow basin}\\{\color{cyan}subspace DACR}}};

\node[shape=ellipse, draw=none, line width=0]() at ({11.5*cos(180)},{11.5*sin(180)}) {\specialcell{{\color{magenta}Neighborhood \&}\\{\color{magenta}Almost Cylinder}\\{\color{cyan}Wide basin}\\{\color{cyan} full space DACR}}};

\node[shape=ellipse, draw=none, line width=0]() at ({11.5*cos(157.5)},{11.5*sin(157.5)}) {\specialcell{{\color{magenta}Cylinder DACR}\\{\color{cyan}Wide basin}\\{\color{cyan}subspace DACR}}};

\node[shape=ellipse, draw=none, line width=0]() at ({10.5*cos(135)},{10.5*sin(135)}) {\color{cyan} \specialcell{Full basin\\DACR}};

\node[shape=ellipse, draw=none, line width=0]() at ({10.5*cos(112.5)},{10.5*sin(112.5)}) {\color{cyan} \specialcell{Full basin\\DACR}};

\node[shape=ellipse, draw=none, line width=0]() at ({10.5*cos(67.5)},{10.5*sin(67.5)}) {\color{cyan} \specialcell{Full basin\\DACR}};

\node[shape=ellipse, draw=none, line width=0]() at ({10.5*cos(45)},{10.5*sin(45)}) {\color{cyan} \specialcell{Full basin\\DACR}};

\node[shape=ellipse, draw=none, line width=0]() at ({10.5*cos(-90)},{10.5*sin(-90)}) {\color{olive} \specialcell{Null DACR}};

\node[shape=ellipse, draw=none, line width=0]() at ({11*cos(-22.5)},{11*sin(-22.5)}) {\color{olive} \specialcell{Null DACR}};

\node[shape=ellipse, draw=none, line width=0]() at ({10.5*cos(202.5)},{10.5*sin(202.5)}) {\color{olive} \specialcell{Null DACR}};

\node[shape=ellipse, draw=none, line width=0]() at ({11*cos(-45)},{11*sin(-45)}) {\color{olive} \specialcell{Null DACR}};

\node[shape=ellipse, draw=none, line width=0]() at ({10.5*cos(-67.5)},{10.5*sin(-67.5)}) {\color{olive} \specialcell{Null DACR}};

\node[shape=ellipse, draw=none, line width=0]() at ({10.5*cos(-112.5)},{10.5*sin(-112.5)}) {\color{olive} \specialcell{Null DACR}};

\node[shape=ellipse, draw=none, line width=0]() at ({11*cos(-135)},{11*sin(-135)}) {\color{olive} \specialcell{Null DACR}};

\node[shape=ellipse, draw=none, line width=0]() at ({3.5*cos(90)},{3.5*sin(90)}) {\color{teal} \specialcell{\Large Network motifs \\\Large with weakly stable\\\Large hyperplane $\{x = x^*\}$}};

\end{tikzpicture}
\caption{{\bf (Motifs of Weak Dynamic ACR)} 
There are 17 motifs of weak dynamic ACR with two reactions and two or fewer species.
A necessary and sufficient condition for weak dynamic ACR (as well as weak full basin dynamic ACR) is that the reactant polytope be parallel to the axis of the ACR variable ({\color{green} green} line segment) and both reactions point inwards. 
Of these motifs, 16 are placed on the circumference of a circle with coordinates $\theta = n\pi/8, n \in \{0,1,\ldots 15\}$, while 1 motif is placed at the center of the circle. The two arrows make the same angle with the reactant polytope for the motifs at $\theta = n \pi/2, n \in \{0,1,2,3\}$ (the four cardinal directions -- north, south, east and west). 
For $\theta \in n\pi/8, n \in \{0,1,2, 3,4\}$ (north-east quadrant of the picture), the left arrow is fixed in the north-east quadrant while the right arrow rotates southwards moving southwards along the picture. Similarly in the north-west quadrant, the right arrow is fixed in the north-west quadrant; in the south-west quadrant the left arrow is fixed in the south-east quadrant; and in the south-east quadrant the right arrow is fixed in the south-west quadrant. 
The figure of motifs is invariant under reflection and under rotation around the central vertical axis. 
The figure is also invariant under a combination of reflection around a central horizontal axis and rotation of each motif around the axis of that motif. 
The central horizontal band -- with the motifs at $\theta = 0, \pi$ on the circumference and the motif at the center of the circle -- have $\dim(\SS)=1$ while the rest have $\dim(\SS)=2$. 
The motif at the center of the circle can have an embedding with either 1 or 2 species (the second species remaining dynamically unchanged), an embedding of every motif on the circumference requires 2 species. 
Each motif is labeled with its strongest local ACR property (in {\color{magenta} magenta}) as well as its strongest non-local ACR property (in {\color{cyan} cyan}). Null ACR is labeled in (in {\color{olive} olive}). Moving northwards along the circumference of the circle, the motifs have stronger local and non-local ACR properties. 
}
\label{fig:p318tyhelgheg}
\end{figure}
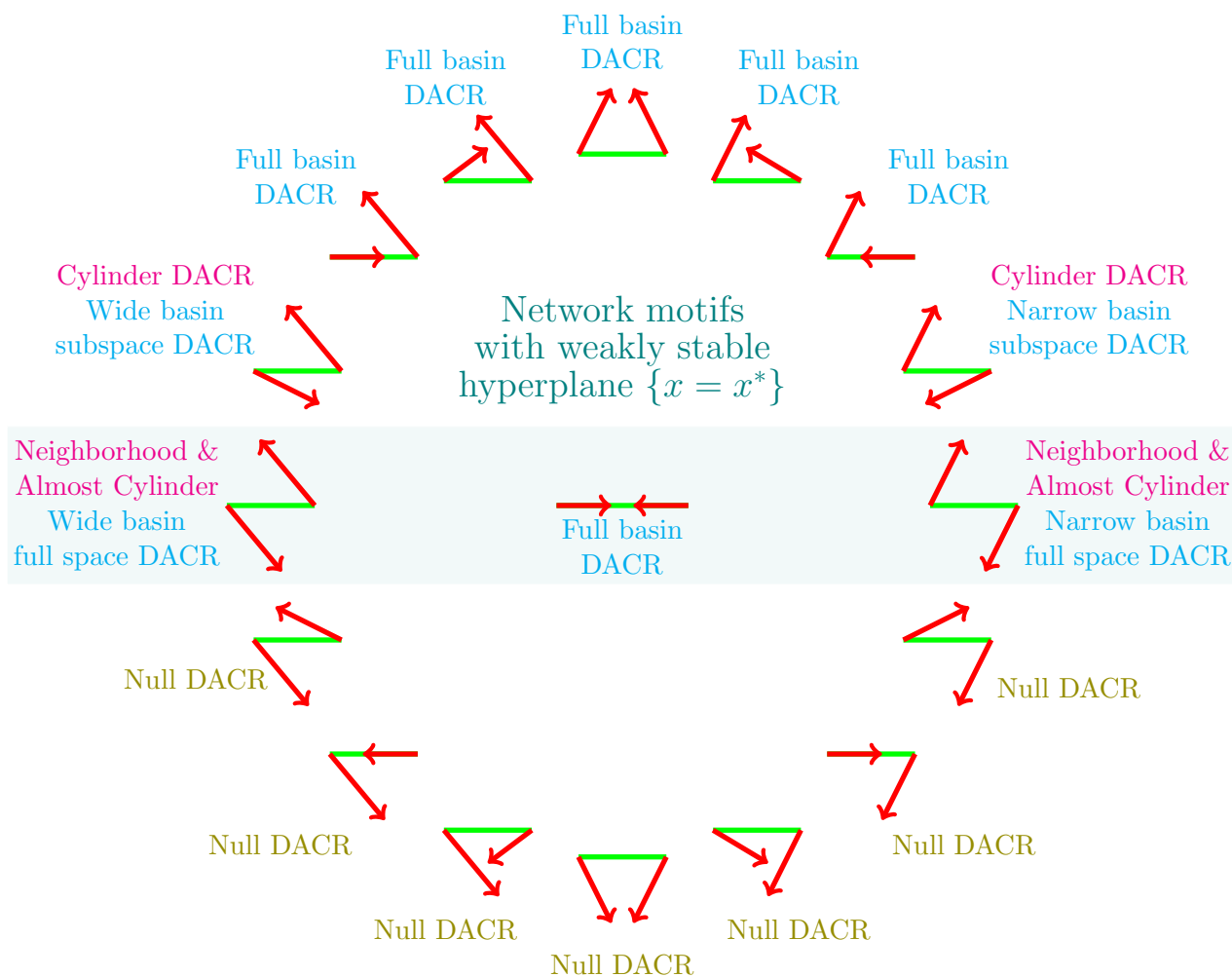

Each network motif is related to an infinite family of networks, see Section \ref{sec:networkmotifs} for an explanation of how a network maps to a motif.
The motifs in Figure \ref{fig:allthemotifs} (or rather the networks that map to the motifs) have a rich diversity of dynamical properties. Key results in this paper, appearing in Section \ref{sec:6yyj56hdgkwtk} onwards, deal with identifying such properties. The results are summarized in Figures \ref{fig:o45yuohdaguh3iur} and \ref{fig:p318tyhelgheg} -- fully annotated versions of Figure \ref{fig:allthemotifs}. The reader is encouraged to begin with at least a glance at {\bf the summary theorem (Theorem \ref{thm:mainesttheorem}) and Figures \ref{fig:o45yuohdaguh3iur} and \ref{fig:p318tyhelgheg}} to orient themselves towards the objectives of the paper.
{\bf Specifically, we identify which of the 8 network motifs with static ACR and which of the 17 network motifs with weak dynamic ACR are also dynamic ACR.} 



Infinitely many reaction networks map onto a single network motif -- an object that captures two crucial network properties common to the family of networks: (i) orientation of the reactant polytope, (ii) orientation of the reaction vectors relative to the coordinate axes and relative to one another. 
These two network properties are sufficient for characterizing a wide range of dynamic ACR properties. 
To understand these network properties, we embed a reaction network in Euclidean space by identifying the stoichiometric coefficients of each species in a complex with a point in Euclidean space. For instance, the reactant and product complex of the reaction $X+Y \to 2Z$ are identified with the points $(1,1,0)$ and $(0,0,2)$ in Euclidean space. 
The reactant polytope of a reaction network is the convex hull of the reactant complexes in Euclidean space. In the case of a reaction network with only two reactions, a reactant polytope is merely a line segment or a single point. 
Reactions are fixed ``arrows'' in Euclidean space. 
For the case of two species, a reaction may be aligned with the coordinate axes or may be in the interior of a quadrant. 
Moreover when there are two reactions, multiple relative orientations of the reaction arrows are possible. 
A network motif is depicted by drawing the orientation of the reactant polytope and the relative orientations of the reactions. 
{\bf A network motif is deeply connected with the dynamics of the reaction system.} For one-dimensional systems such as the one in Figure \ref{fig:qperuhg305y84oigg}(a), the reaction vectors must be parallel to one another and moreover the trajectories in phase space, when superimposed on the reaction network embedding, are simply parallel to the reaction vectors. For two and higher dimensional systems, the relationship is more complicated since the vector field is a linear combination of the reaction vectors with variable coefficients that depend on the position in phase space. 
The reactant polytope is extremely significant for ACR. 
All the various forms of ACR properties with at most two species require that the reactant polytope be parallel to one of the coordinate axes, i.e. the polytope should be a horizontal or a vertical line segment. When the reactant polytope is horizontal, the ACR hyperplane is vertical and vice versa. 
Among static ACR networks in two species, we show that dynamic ACR requires that the reactions points ``inwards'', along the reactant polytope. 
For the weak dynamic ACR networks in two species, we show that dynamic ACR requires that on average, reactions point away from the coordinate axis in the direction perpendicular to the reactant polytope (``up'' when the reactant polytope is horizontal). 
Of course, the above statements are merely suggestive, see Theorem \ref{thm:mainesttheorem} for precise statements. 

\vspace{0.1in}
\noindent{\large \bf Related work:}
The idea of minimal reaction networks with a certain dynamical property has a long history. 
While this work is the first to identify minimal motifs of dynamic ACR, network motifs for other dynamical properties have been studied extensively. 
For example, networks with two species which have limit cycles have been identified in \cite{schnakenberg1979simple}, while \cite{schlogl1972chemical} studies simple  networks with  phase transitions giving candidates for multistationarity. 
In recent years, several classes of minimal multistationary networks called {\em atoms of multistationarity} have been identified in \cite{joshi2013atoms,joshi2013complete,joshi2015survey,craciun2021multistationarity}, and those of {\em oscillations} in \cite{banaji2018inheritance2}. 
Network conditions for static ACR and minimal networks with static ACR have also appeared in recent work \cite{shinar2011design,pascual2020local,meshkat2021absolute}. 
Biochemical implications of static ACR have been studied \cite{karp2012complex,dexter2013dimerization,dexter2015invariants} as well as the implications from a control theory perspective \cite{cappelletti2020hidden}. 
Stochastic (continuous-time Markov chain) models of reaction networks with the ACR property were studied in \cite{anderson2014stochastic,anderson2017finite}. 


This article is organized as follows.   
Section~\ref{sec:formsofacr} introduces different forms of static and dynamic ACR and the relations between them. 
Section~\ref{sec:6yyj56hdgkwtk} gives a classification of small reaction networks, those containing at most two reactions and two species; here the focus is on mass-conserving reaction networks. 
Section~\ref{sec:dynamicacr_classification} continues the classification, extended to networks where the trajectories may go to infinity or to the boundary. 
Section~\ref{sec:summaryofallacr} summarizes the main classification results with the statement of a larger theorem that encapsulates numerous results from the previous sections.

\section{Forms of ACR} \label{sec:formsofacr}

\subsection*{Background and terminology}
We use standard notation and terminology for reaction networks and mass action systems. 
A quick summary is given here, see for instance \cite{joshi2015survey,yu2018mathematical} for further details. 
Upper case letters ($X,Y, Z, A, B$) are used for species participating in reactions and the corresponding lower case letters ($x,y,z,a,b$) for their concentrations which are time-varying quantities. An example of a reaction is $X + Y \to 2Z$, where $X+Y$ is referred to as the source complex, while $2Z$ is the product complex. The rate of any given reaction is a nonnegative-valued function of species concentrations. In the case of mass action kinetics, the rate is proportional to the product of reactant concentrations taken with multiplicity. The proportionality constant, called the reaction rate constant, is placed adjacent to the reaction arrow, as follows: $X + Y \xrightarrow{k} 2Z$. The rate of this reaction under mass action kinetics is $kxy$. 
The reaction vector for this reaction is the difference between the product complex and the source complex, i.e. $2Z-(X+Y)$, which under a choice of a standard coordinate basis is written as $(-1,-1,2)$. 
A reaction network is a nonempty set of reactions, such that every species participates in at least one reaction, and none of the reaction vectors is the zero vector. 
The {\em stoichiometric subspace} of a reaction network is the subspace spanned by the set of reaction vectors of the reaction network.


We use $\GG$ to denote a reaction network and $K$ to denote a specific choice of mass action kinetics for $\GG$, so that $(\GG,K)$ is a mass action dynamical system. 
Throughout the paper, we depict the Euclidean embedding of a reaction network next to the dynamics in the phase plane. Figure \ref{fig:qperuhg305y84oigg}(a) is an example of a Euclidean embedding: the reactions $B \to A$ and $A+B \to 2B$ are shown as red arrows, each arrow originating at the reactant complex and terminating at the product complex. We also depict the {\em reactant polytope}, defined as the convex hull of the reactant complexes. In the example in Figure \ref{fig:qperuhg305y84oigg}(a), the set of reactant complexes is $\{B, A+B\}$ and so the reactant polytope, shown in green, is a line segment joining the two complexes.

Throughout this paper, we consider a dynamical system $\DD$ defined by $\dot x = f(x)$ with $x \in \R^n_{\ge 0}$ for which $\R^n_{\ge 0}$ is forward invariant. $x \in \R^n_{\ge 0}$ is a {\em steady state} of $\DD$ if $f(x)=0$. 

\begin{definition} \label{def:0496yuwpeh}
The {\em kinetic subspace} of $\DD$ is defined to be the linear span of the image of $f$, denoted by $\spn(\Im(f))$. 
$x,y \in \R^n_{\ge 0}$ are {\em compatible} if $y - x \in \spn(\Im(f))$. 
The sets $S, S' \subseteq \R^n_{\ge 0}$ are {\em compatible} if there are $x \in S$ and $x' \in S'$ such that $x$ and $x'$ are compatible.
\end{definition}


The notation $\HH[i, a_i^*] \coloneqq \{x \in \R^n_{> 0}: x_i = a_i^*\}$ is reserved for the hyperplane parallel to a coordinate hyperplane and restricted to the positive orthant.  
When variables are labelled without indices, for instance as $x, y, z$ etc., we use the notation $\HH[x, x^*] \coloneqq \{x \in \R^n_{> 0}: x = x^*\}$. 
Even though the two notations are slightly inconsistent, there is no possibility of confusion.

Static ACR and dynamic ACR in a real dynamical system were defined in \cite{joshi2021foundations}. 
We repeat these definitions here: 
\begin{definition} \label{def:42u5hgd;;goi}
\beit
\item $\DD$ is a {\em static ACR system} if $\DD$ has a positive steady state and there is an $i \in \{1,\ldots, n\}$ and a positive $a_i^* \in \R_{> 0}$ such that any positive steady state $x \in \R^n_{> 0}$ is in the hyperplane $\{x_i = a_i^*\}$. 
Any such $x_i$ and $a_i^*$ is a {\em static ACR variable} and its {\em static ACR value}, respectively. $\HH[i,a_i^*]$ is the static ACR hyperplane. 
\item $\DD$ is a {\em dynamic ACR system} if there is an $i \in \{1,\ldots, n\}$ with $f_i \not \equiv 0$ and a positive $a_i^* \in \R_{> 0}$ such that for any $x(0) \in \R^n_{> 0}$ that is compatible with $\{x \in \R^n_{>0} ~|~ x_i = a_i^*\}$, a unique solution to $\dot x = f(x)$ exists up to some maximal $T_0(x(0)) \in (0, \infty]$, and $x_i(t) \xrightarrow{t \to T_0} a_i^*$. Any such $x_i$ and $a_i^*$  is a {\em dynamic ACR variable} and its {\em dynamic ACR value}, respectively.  $\HH[i,a_i^*]$ is the dynamic ACR hyperplane.
\enit
\end{definition}
Now we define strong and weak forms of each type of ACR. 
All forms: static ACR, strong static ACR, dynamic ACR and weak dynamic ACR are marked by the presence of a privileged hyperplane parallel to a coordinate hyperplane $\{x \in \R^n_{> 0}: x_i = a_i^*\}$ called the {\bf ACR hyperplane}. 
This hyperplane either contains all the steady states (in the static case) or is the unique attractor for all relevant trajectories (in the dynamic case). 
Moreover, each form of ACR has a set $\Omega$ associated with it called the {\bf basin of ACR}. 

\subsection{Static forms of ACR}

We generalize the definition of static ACR from \cite{joshi2021foundations} to allow for arbitrary basin sets $\Omega$. 
\begin{definition} \label{def:qe;orghqeorghokg}
The variable $x_i$, where $i \in \{1, \ldots, n\}$, has {\em static ACR w.r.t. $\Omega \subseteq  \R_{\ge 0}^n$}  if there is an $a_i^* > 0$ such that the following hold:
\been
\item $f(x) = 0$ for some positive $x \in \Omega$,  
\item for any $x \in \Omega$ such that $f(x)=0$, $x_i = a_i^*$. 
\enen
In this case, the {\em static ACR value} of $x_i$ is $a_i^*$ and $\HH[i, a_i^*] \coloneqq \{x \in \R^n_{> 0}: x_i = a_i^*\}$ is the {\em static ACR hyperplane}. 
\end{definition}

\begin{remark}
When $\Omega = \R^n_{> 0}$, we simply say that a variable has ``static ACR'' instead of ``static ACR w.r.t. $\R^n_{> 0}$''. Note that boundary steady states are conventionally excluded from consideration, in other words, we do not require that $\Omega = \R^n_{\ge 0}$ for static ACR, merely that $\Omega = \R^n_{> 0}$. 
\end{remark}

By strong static ACR, we mean the property that every compatibility class that intersects the static ACR hyperplane contains at least one positive steady state. 
Most commonly studied motifs and biochemical systems with static ACR are strong static ACR as well.

\begin{definition} \label{def:djlu5feiulfhktyu}
The variable $x_i$, where $i \in \{1, \ldots, n\}$, has {\em strong static ACR w.r.t. $\Omega \subseteq \R_{\ge 0}^n$}  if there is an $a_i^* > 0$ such that the following hold:
\been
\item $x_i$ has static ACR w.r.t. $\Omega \subseteq \R_{\ge 0}^n$ with value $a_i^*$,  
\item for any $y \in \HH[i, a_i^*] \cap \Omega$  there is a $z \in \Omega \cap \R^n_{> 0}$ such that $y-z \in \SS_f$ and $f(z) = 0$.
\enen
In this case, the {\em strong static ACR value} of $x_i$ is $a_i^*$ and $\HH[i, a_i^*]$ is the {\em strong static ACR hyperplane}. 
\end{definition}

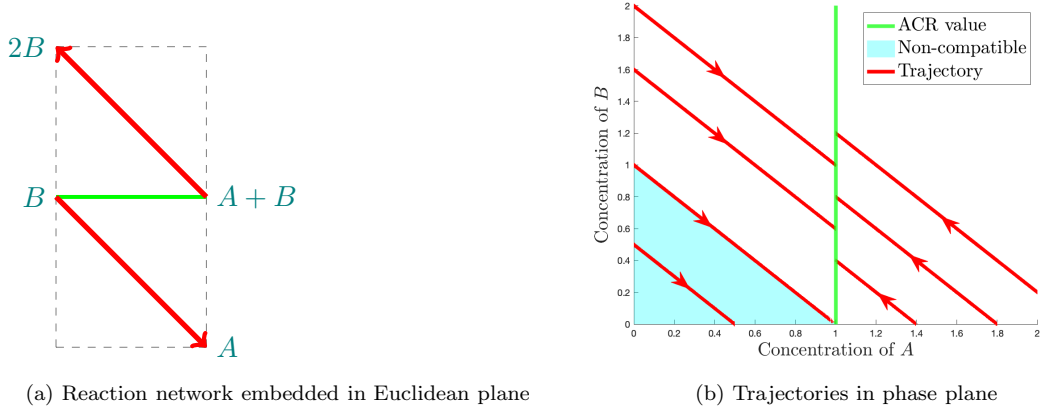
\begin{figure}[h!] 
\centering
\begin{subfigure}[b]{0.45\textwidth}
\begin{tikzpicture}[scale=2]
\draw[help lines, dashed, line width=0.25] (0,0) grid (1,2);

\node [right] at (1,0) {{\color{teal} $A$}};
\node [left] at (0,1) {{\color{teal} $B$}};
\node [right] at (1,1) {{\color{teal} $A+B$}};
\node [left] at (0,2) {{\color{teal} $2B$}};

\draw [-, line width=1.5, green] (0,1) -- (1,1);
\draw [->, line width=2, red] (0,1) -- (1,0);
\draw [->, line width=2, red] (1,1) -- (0,2);

\end{tikzpicture}
\caption{Reaction network embedded in Euclidean plane}
\end{subfigure}
\begin{subfigure}[b]{0.45\textwidth}
\includegraphics[scale=0.35]{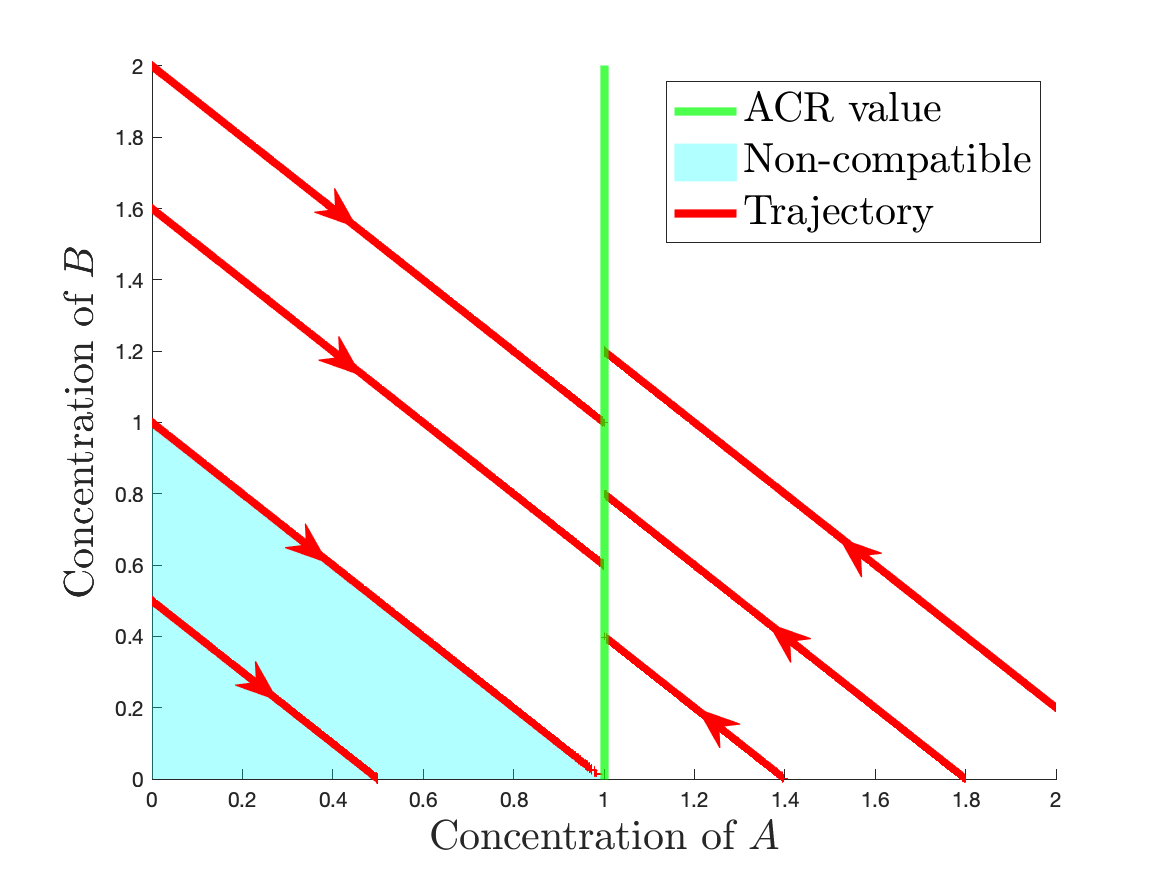}
\caption{Trajectories in phase plane}
\end{subfigure}
\caption{{\bf (Archetypal wide basin dynamic ACR network)} A dynamic ACR reaction network ($A+B \to 2B, B \to A$) with $A$ as a wide basin dynamic ACR variable. 
The concentration of $A$ is bounded within the subset of $\R^2_{\ge 0}$ that is not compatible the ACR hyperplane $\{a=1\}$ (non-compatible region shown here in cyan).}
\label{fig:qperuhg305y84oigg}
\end{figure}

The concentration of $A$ in the network in Figure \ref{fig:qperuhg305y84oigg} has strong static ACR w.r.t. $\R^2_{> 0}$ because every point on the hyperplane $\HH[a,a^*] = \{(a,b) \in \R^2_{> 0} : a=a^*\}$ (the green vertical line) is a steady state. So any compatibility class with $a(0) + b(0) > a^*$ intersects the hyperplane $\HH[a,a^*]$ and therefore has a positive steady state. 
Every network motif with static ACR studied in this paper also has strong static ACR, as we will show in Theorem \ref{thm:po43ghi3lwkfh}.


\subsection{Dynamic forms of ACR}
We generalize the definition of dynamic ACR from \cite{joshi2021foundations} in two ways: (i) we allow for arbitrary basin sets $\Omega$, and (ii) we define a weaker form for which the ACR hyperplane is weakly attracting. 
\begin{definition} \label{def:fjhdkurytjkayjghm,kuzf.n}
The variable $x_i$, where $i \in \{1, \ldots, n\}$, has {\em dynamic ACR w.r.t. $\Omega \subseteq \R_{\ge 0}^n$} if there is an $a_i^* > 0$ such that for any initial value $x(0)$ in $\Omega$, a unique solution to $\dot x = f(x)$ exists up to some maximal $T_0(x(0)) \in (0, \infty]$, and $\lim_{t \to T_0} x_i(t) = a_i^*$. 
We say that $x_i$ has {\em dynamic ACR value} $a_i^*$. Moreover, we say that the ACR hyperplane $\HH[i,a_i^*]$ is an {\em attractor} for $\Omega$ and that $\Omega$ is a {\em basin of attraction} of  $\HH[i,a_i^*]$. 
\end{definition}


Weak dynamic ACR is the notion that the ACR variable converges to a value that is not further from the ACR value than the initial distance. In other words, all initial conditions are attracted towards the ACR hyperplane even if they fail to reach there. 
\begin{definition} \label{def:w4ophujethlndlnl}
The variable $x_i$, where $i \in \{1, \ldots, n\}$, has {\em weak dynamic ACR w.r.t. $\Omega \subseteq \R^n_{\ge 0}$} if there is an $a_i^* > 0$ such that for any initial value $x(0)$ in $\Omega$, a unique solution to $\dot x = f(x)$ exists up to some maximal $T_0(x(0)) \in (0, \infty]$,  
$\lim_{t \to T_0} x_i(t)$ exists and 
$\lim_{t \to T_0}\abs{x_i(t) - a_i^*} \le \abs{x_i(0) - a_i^*}$, with strict inequality when $x_i(0) \ne a_i^*$. 
We say that $x_i$ has {\em weak dynamic ACR value} $a_i^*$. Moreover, we say that the ACR hyperplane $\HH[i,a_i^*]$ is a {\em weak attractor} for $\Omega$ and $\Omega$ is a {\em weak basin of attraction} of  $\HH[i,a_i^*]$. 
\end{definition}

\begin{remark}
We do not require the (weak) ACR hyperplane is invariant, merely that it is a (weak) attractor. A weak attractor may fail to be an attractor because trajectories reach the boundary or diverge to infinity before reaching the ACR hyperplane.  
\end{remark}

%

In analogy with the definition of stability of a steady state, we define the following.

\begin{definition} \label{def:o4tyhoegjwppp} 
\beit
\item
The hyperplane $\HH[i,a_i^*]$ is {\em stable} w.r.t. $\Omega$ if $\HH[i,a_i^*]$ is an attractor for $\Omega$ and all trajectories with initial value in $\Omega$ move monotonically towards $\HH[i,a_i^*]$. 
\item
The hyperplane $\HH[i,a_i^*]$ is {\em weakly stable} w.r.t. $\Omega$ if $\HH[i,a_i^*]$ is a weak attractor for $\Omega$ and all trajectories with initial value in $\Omega$ move monotonically towards $\HH[i,a_i^*]$. 
\enit
\end{definition}

\begin{remark}
Stability (weak or not) of $\HH[i,a_i^*]$ implies that $\HH[i,a_i^*]$ is invariant. 
\end{remark}

\begin{example}{\em (Weakly stable but not stable hyperplane)}
Consider the mass action system of the following reaction network (see also Figure \ref{fig:o;eithjeoihy3j2potj}(a)): 
\alis{
2A+B \xrightarrow{k_1} 2B, \quad B \xrightarrow{k_2} A. 
}
The mass action system of ODEs is:
\begin{align} \label{eq:;q3hio3hirgo3j}
\dot a = b(k_2 -2k_1 a^2), \quad \dot b = -b(k_2 - k_1 a^2). 
\end{align}
It is clear from $\dot a$ that for the duration of time that $b(t)$ is positive, $a(t)$ approaches the value $\sqrt{k_2/(2k_1)}$.
This implies that $a$ is a weak dynamic ACR variable with a weak ACR value of $\sqrt{k_2/(2k_1)}$. Note however that $\dot a + \dot b = -k_1a^2b$ is negative everywhere in the positive orthant, and so we expect $b(t)$ to converge to $0$ for most initial values. In fact, this is the case for every initial value that is not on the weak ACR hyperplane, as a consequence of Theorem \ref{thm:qotih3ohoeghow}.  
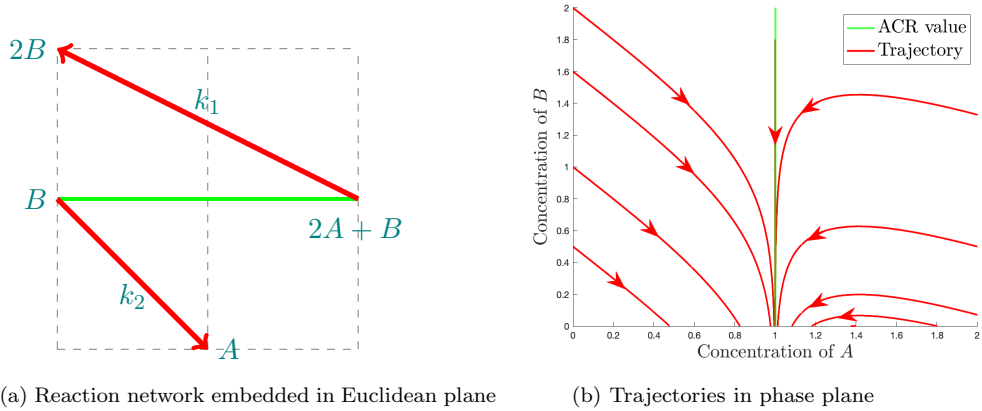
\begin{figure}[h!] 
\centering
\begin{subfigure}[b]{0.4\textwidth}
\begin{tikzpicture}[scale=2]
\draw[help lines, dashed, line width=0.25] (0,0) grid (2,2);

\node [right] at (1,0) {{\color{teal} $A$}};
\node [left] at (0,1) {{\color{teal} $B$}};
\node [right] at (1.6,0.8) {{\color{teal} $2A+B$}};
\node [left] at (0,2) {{\color{teal} $2B$}};

\node [above] at (1,1.5) {{\color{teal} $k_1$}};
\node [below] at (.5,0.5) {{\color{teal} $k_2$}};

\draw [-, line width=1.5, green] (0,1) -- (2,1);
\draw [->, line width=2, red] (0,1) -- (1,0);
\draw [->, line width=2, red] (2,1) -- (0,2);

\end{tikzpicture}
\caption{Reaction network embedded in Euclidean plane}
\end{subfigure}
\begin{subfigure}[b]{0.35\textwidth}
\includegraphics[scale=0.35]{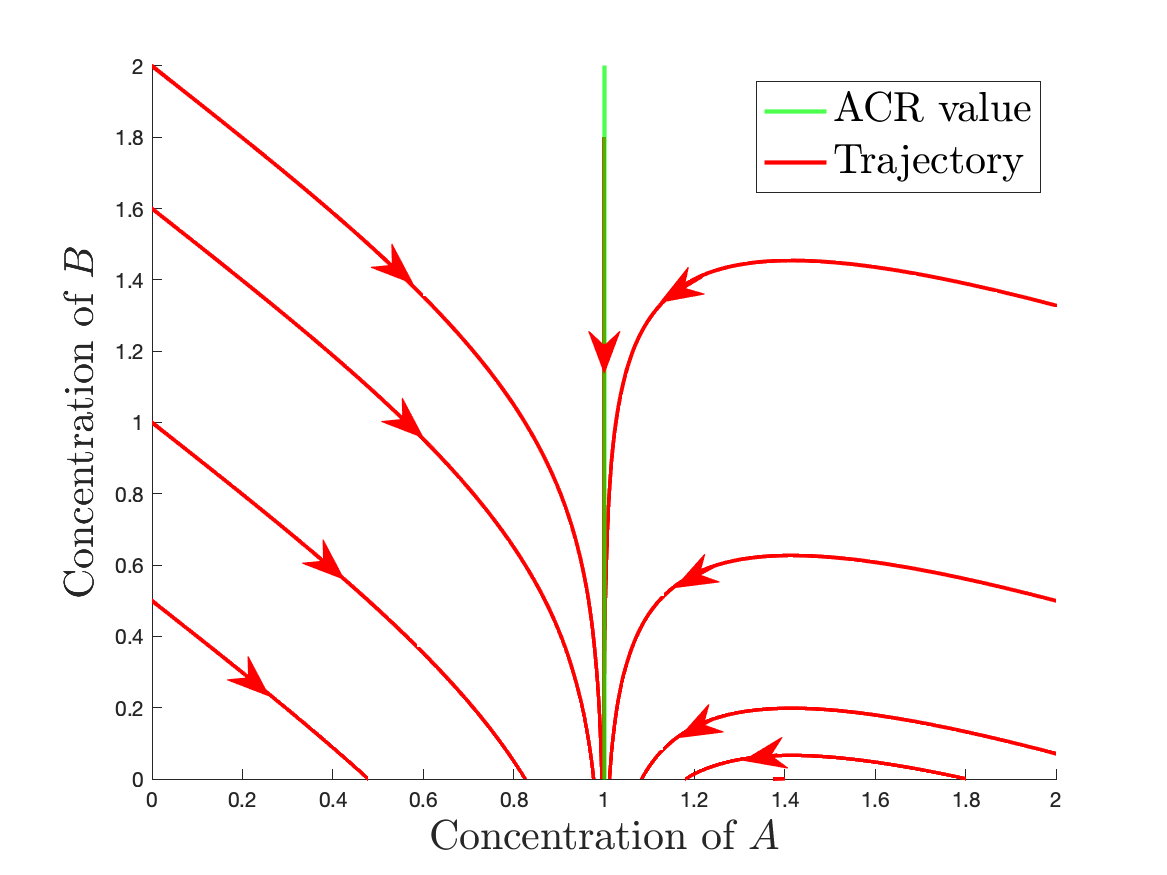}
\caption{Trajectories in phase plane}
\end{subfigure}
\caption{The mass action system $\{2A+B \xrightarrow{k_1} 2B,  B \xrightarrow{k_2} A\}$, for any choice of rate constants, has only weak dynamic ACR in the concentration of $A$. All trajectories move towards the ACR hyperplane, but do not approach the ACR hyperplane in the limit of large time.}
\label{fig:o;eithjeoihy3j2potj}
\end{figure}


\end{example}

\subsection{Relations between different forms of ACR}

It is clear from the definitions that if $x_i$ has dynamic ACR w.r.t. $\Omega$ with value $a_i^*$, then $x_i$ has weak dynamic ACR w.r.t. $\Omega$ with the same value $a_i^*$. Similarly strong static ACR w.r.t. $\Omega$ implies static ACR w.r.t. $\Omega$. The static forms of ACR are also related to the dynamic forms under some mild assumptions of existence of steady states. 

\begin{theorem} \label{thm:formsofacr}
Suppose that $x_i$ is weak dynamic ACR w.r.t. $\Omega$ with value $a_i^*$. 
\been
\item Suppose there is an $x \in \Omega$ such that $f(x) =0$. Then $x_i$ is static ACR w.r.t. $\Omega$.
\item Suppose for any $y \in \HH[i, a_i^*] \cap \Omega$  there is a positive $x \in \Omega$ such that $y-x \in \SS_f$ and $f(x) = 0$. Then $x_i$ is strong static ACR w.r.t. $\Omega$.
\enen
\end{theorem}
\begin{proof}
For both cases, the static ACR property follows from observing that there cannot be any positive steady state in $\Omega \setminus \HH[i,a_i^*]$ because such a steady state violates the weak dynamic ACR hypothesis. The additional implication of strong static ACR in the second case is immediate from the definition. 
\end{proof}

\begin{figure}[h!] 
 \begin{center}
   \begin{tikzpicture}[auto, every node/.style={scale=0.9}]
        \tikzstyle{block} = [draw, rectangle];
        \node [block] (dyn) {\specialcell{$\{x=x^*\}$ attracts $\Omega$\\($x$ dynamic ACR w.r.t. $\Omega$)}};
         \node [block, left=2cm of dyn] (stab) {$\{x=x^*\}$ stable w.r.t. $\Omega$};
        \node [block, below=1.1cm of dyn] (weak) {\specialcell{$\{x=x^*\}$ weakly attracts $\Omega$\\($x$ weak dynamic ACR w.r.t. $\Omega$)}};
        \node [block, below=1.5cm of stab] (weakstab) {\specialcell{$\{x=x^*\}$ weakly stable w.r.t. $\Omega$}};
        \node [block, below right=1.0cm of weak] (strongstatic) {$x$ strong static ACR w.r.t. $\Omega$};
        \node [block, below=1.25cm of strongstatic] (static) {$x$ static ACR w.r.t. $\Omega$};

     \draw[-implies,double equal sign distance] (dyn)--(weak);
     \draw[-implies,double equal sign distance] (strongstatic)--(static);
     \draw[-implies,double equal sign distance] (stab)--(weakstab);
     \draw[-implies,double equal sign distance] (stab)--(dyn);
     \draw[-implies,double equal sign distance] (weakstab)--(weak);
          {\color{cyan}
     \draw[-implies,double equal sign distance] (weak)--(static);
     }
     {\color{magenta}
     \draw[-implies,double equal sign distance] (weak)--(strongstatic);
     }

   \end{tikzpicture} 
 \end{center}
 \caption{({\bf Relations between static ACR and dynamic ACR}) The implications in black follow from the definitions. The implication in {\color{cyan} cyan} requires the additional assumption of existence of a positive steady state. The implication in {\color{magenta} magenta} requires the additional assumption of existence of a positive steady state within each compatibility class which has a nonempty intersection with the ACR hyperplane $[i, a_i^*]$.}
 \label{fig:3oyho3thwef;w;o}
 \end{figure}
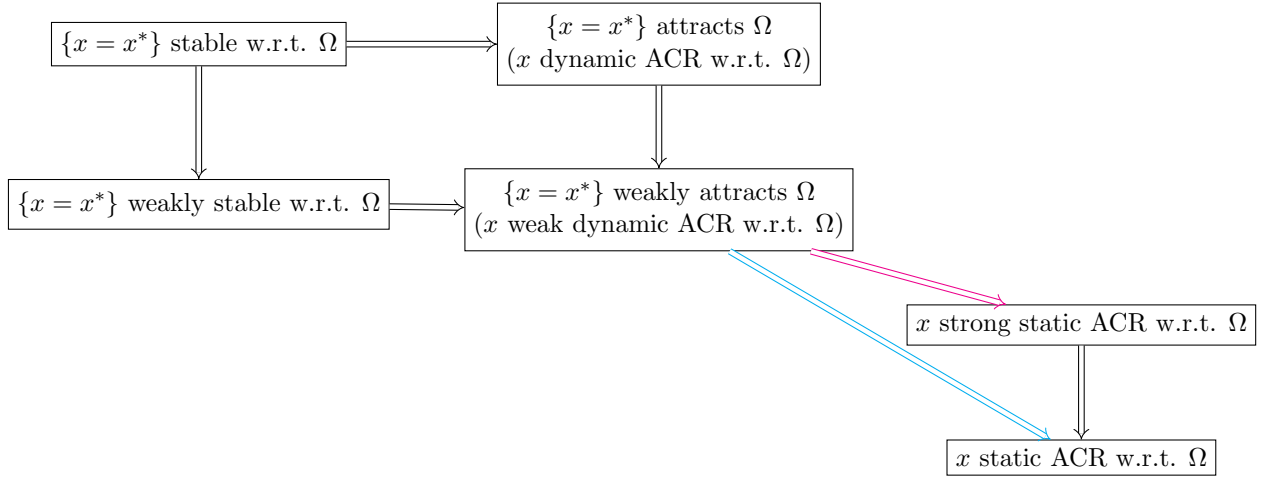

These relations are portrayed in Figure \ref{fig:3oyho3thwef;w;o}. 


\subsection{Some basins of interest for all forms of ACR}

We discuss some basins of natural interest and the relations between them. A basin of ACR applies to any of the forms of ACR (static, strong static, dynamic, weak dynamic) discussed earlier.  

\begin{definition}
Let $\PP \in \{$static, strong static, dynamic, weak dynamic$\}$. We define various {\em basin types} as follows. 
\begin{enumerate}[(i)]
\item {\em Full basin $\PP$-ACR} occurs when $\Omega = \R^n_{> 0}$. {\em Full basin static ACR} is simply referred to as {\em static ACR}. 

\item {\em Subspace $\PP$-ACR} occurs when $\Omega = (\HH[i,a_i^*] + \SS) \cap \R^n_{> 0}$ for some subspace $\SS$ of $\SS_f$ such that $\SS \not \subseteq e_i^\perp$.  
We say {\em full space $\PP$-ACR} if we have subspace $\PP$-ACR with $\SS = \SS_f$. 
Full space dynamic ACR is simply referred to as {\em dynamic ACR}. 

\item {\em Neighborhood $\PP$-ACR} occurs when $\Omega$ is a neighborhood $\HH[i,a_i^*]$. 

Suppose there are some $M_j > 0$ for all $j \ne i$ such that the neighborhood of the set $\{x_i = a_i^*, x_j > M_j: j \ne i\}$ is a basin of $\PP$-ACR. Then we say that $x_i$ has {\em almost neighborhood $\PP$-ACR}.


\item {\em Cylinder $\PP$-ACR} occurs when $\Omega$ is a cylinder of $\HH[i,a_i^*]$. A cylinder of $\HH[i,a_i^*]$ is the set of points $\{\abs{x_i - a_i^*} < \delta^*\}$ for some $\delta^*>0$. 

We define {\em almost cylinder $\PP$-ACR} when a basin of $\PP$-ACR is a cylinder of some set $\{x_i = a_i^*, x_j > M_j: j \ne i\}$ for some $M_j > 0$ for all $j \ne i$. 


\item {\em Null $\PP$-ACR} occurs when $\Omega = \HH[i,a_i^*]$. 
{\em Non-null $\PP$-ACR} occurs when $\Omega \setminus \HH[i,a_i^*] \ne \varnothing$.  

\enen
\end{definition}

 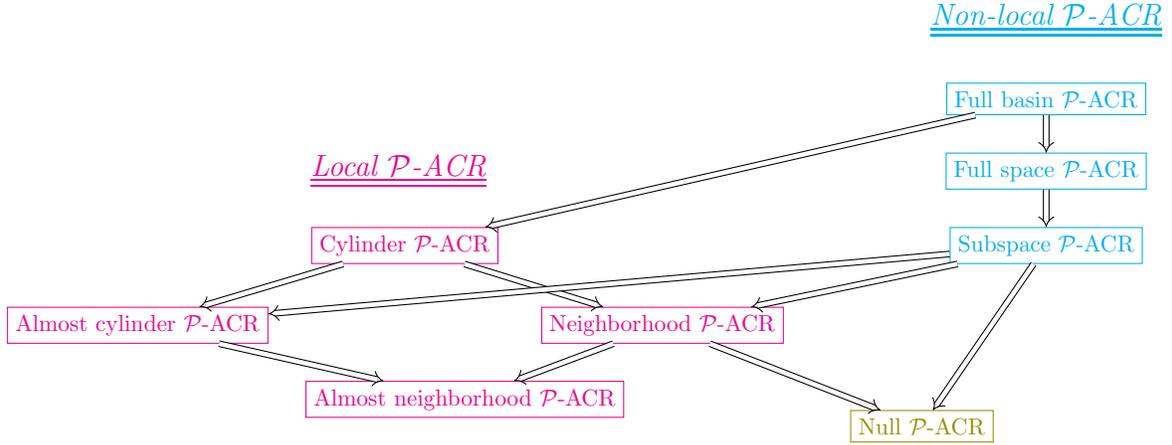
\begin{figure}[h!]
 \begin{center}
   \begin{tikzpicture}[auto, every node/.style={scale=0.8}]
        \tikzstyle{block} = [draw, rectangle];
        {\color{cyan}
        \node [] (nonlocal) {\Large \em \underline{\underline{Non-local $\PP$-ACR}}};
        \node [block, below=0.5cm of nonlocal] (fullbasin) {Full basin $\PP$-ACR};
        \node [block, below=0.5cm of fullbasin] (dynamic) {Full space $\PP$-ACR};
        \node [block, below=0.5cm of dynamic] (subspace) {Subspace $\PP$-ACR};
        }

        {\color{magenta}
        \node [left=6cm of dynamic] (local) {\Large \em \underline{\underline{Local $\PP$-ACR}}};
         \node [block, left=6cm of subspace] (cylinder) {Cylinder $\PP$-ACR};
         \node [block, below right=.807cm of cylinder] (neighborhood) {Neighborhood $\PP$-ACR};
         \node [block, below left=.807cm of cylinder] (almostcylinder) {Almost cylinder $\PP$-ACR};
         \node [block, below right=.7cm of almostcylinder] (almostneighborhood) {Almost neighborhood $\PP$-ACR};
         }
         
        {\color{olive}
        \node [block, below right=1.25cm of neighborhood] (null) {Null $\PP$-ACR};
        }
        
     \draw[-implies,double equal sign distance] (fullbasin)--(dynamic);
     \draw[-implies,double equal sign distance] (dynamic)--(subspace);
     \draw[-implies,double equal sign distance] (subspace)--(null);
     \draw[-implies,double equal sign distance] (cylinder)--(neighborhood);
      \draw[-implies,double equal sign distance] (cylinder)--(almostcylinder);
       \draw[-implies,double equal sign distance] (neighborhood)--(almostneighborhood);
        \draw[-implies,double equal sign distance] (almostcylinder)--(almostneighborhood);
        
            \draw[-implies,double equal sign distance] (subspace)--(neighborhood);
               \draw[-implies,double equal sign distance] (subspace)--(almostcylinder);
                  \draw[-implies,double equal sign distance] (fullbasin)--(cylinder);
                    \draw[-implies,double equal sign distance] (neighborhood)--(null);
   \end{tikzpicture} 
 \end{center}
 \caption{({\bf Relations between basin types.}) Let $\PP \in \{$static, strong static, dynamic, weak dynamic$\}$. The basin type implications are based on the observation that if $\Omega \subseteq \Omega'$ then $\PP$-ACR w.r.t. $\Omega'$ implies $\PP$-ACR w.r.t. $\Omega$.}
 \label{fig:basintypes}
 \end{figure}
 

Neighborhood ACR, cylinder ACR, almost neighborhood ACR and almost cylinder ACR are {\em local forms of ACR}. Full basin and subspace ACR are {\em non-local forms of ACR}. Null $\PP$-ACR is considered neither local nor non-local. 
Conversely, all local and non-local forms of $\PP$-ACR are non-null. 

It's clear that if $\Omega \subseteq \Omega'$ then $\PP$-ACR w.r.t. $\Omega'$ implies $\PP$-ACR w.r.t. $\Omega$. 
Certain relations between non-local forms: (full basin $\PP$-ACR $\implies$ full space $\PP$-ACR $\implies$ subspace $\PP$-ACR $\implies$ null $\PP$-ACR) and between local forms: (cylinder $\PP$-ACR $\implies$ neighborhood $\PP$-ACR $\implies$ null $\PP$-ACR \& almost neighborhood $\PP$-ACR) and (cylinder $\PP$-ACR $\implies$ almost cylinder $\PP$-ACR $\implies$ almost neighborhood $\PP$-ACR) follow from the definitions. Certain local forms are related to subspace $\PP$-ACR as the following theorem shows. 

\begin{theorem} \label{thm:egopepghoeghoe}
Let $\DD$ be a dynamical system which is subspace $\PP$-ACR. Then $\DD$ is (i) almost cylinder $\PP$-ACR and (ii) neighborhood $\PP$-ACR.
\end{theorem}
\begin{proof}
Let $\SS$ be a subspace of the stoichiometric space $\SS_f$ such that $x_1$ is a subspace $\PP$-ACR variable with value $a_1^*$. 
Let $v \in \SS \cap \left(\R^n \setminus e_1^\perp\right)$, i.e. $v = (v_1,\ldots, v_n)$ is some vector with $v_1 \ne 0$.
Let $\varepsilon \in (0, a_1^*)$ and 
define the almost cylinder neighborhood of $\HH[i,a_i^*]$: 
\[
\Omega_\ve \coloneqq \left \{z \in \R^n_{> 0} : \abs{z_1 - a_1^*} < \varepsilon, ~z_2 > \ve \frac{\abs{v_2}}{\abs{v_1}}, \ldots ,  z_n > \ve \frac{\abs{v_n}}{\abs{v_1}} \right \}.
\] 
To see $x_1$ is $\PP$-ACR w.r.t. $\Omega_\ve$ for every $\ve \in (0, a_1^*)$, we only need to show that $\Omega_\ve$ is contained in $\HH[1,a_1^*] + \spn\{v\}$ which in turn is clearly contained in $\HH[1,a_1^*] + \SS$. 
Indeed, let $z \in \Omega_\ve$. 
Then $z \in \HH[1,a_1^*] + \spn\{v\}$ if there is a $\beta \in \R$ such that $z - \beta v \in \HH[1,a_1^*]$. 
Let $\beta \coloneqq (z_1 - a_1^*)/v_1$. Then 
\[
z - \beta v = \left(a_1^*, z_2 - \frac{z_1-a_1^*}{v_1} v_2, \ldots, z_n - \frac{z_1-a_1^*}{v_1} v_n \right).
\]
For $j \in \{2, \ldots, n \}$, 
\[
\abs{\frac{z_1-a_1^*}{v_1} v_j} = \abs{z_1-a_1^*} \frac{\abs{v_j}}{\abs{v_1}}  < \ve \frac{\abs{v_j}}{\abs{v_1}} < z_j, 
\]
which shows that all components of $z - \beta v$ are positive and the first component is $a_1^*$. 
This proves the almost cylinder property where the almost cylinder is $\Omega_\ve$ for any $\ve \in (0, a_1^*)$. 

Let $\Omega \coloneqq \cup_{\ve > 0} \Omega_\ve$.
It is clear that $\Omega$ is contained in $\HH[1,a_1^*] + \SS$ and that $\Omega$ is a neighborhood of $\HH[1,a_1^*]$. 
This shows that $x_1$ is neighborhood $\PP$-ACR. 
%
%
%
%
\end{proof}
The relations between $\PP$-ACR with different basin types are depicted in Figure \ref{fig:basintypes}. 


\subsection{Some examples of systems with local dynamic ACR}

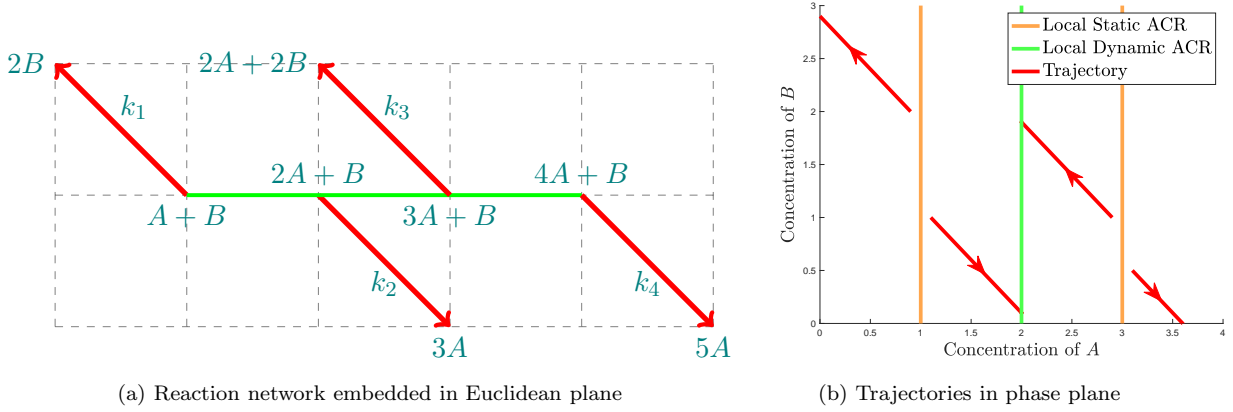
\begin{figure}[h!] 
\centering
\begin{subfigure}[b]{0.6\textwidth}
\begin{tikzpicture}[scale=1.75]
\draw[help lines, dashed, line width=0.25] (0,0) grid (5,2);
\node [below] at (1,1) {{\color{teal} $A+B$}};
\node [left] at (0,2) {{\color{teal} $2B$}};
\node [above] at (2,1) {{\color{teal} $2A+B$}};
\node [below] at (3,0) {{\color{teal} $3A$}};
\node [above] at (.6,1.5) {{\color{teal} $k_1$}};
\node [below] at (2.5,0.5) {{\color{teal} $k_2$}};
\draw [->, line width=2, red] (1,1) -- (0,2);
\draw [->, line width=2, red] (2,1) -- (3,0);
\draw [-, line width=1.5, green] (1,1) -- (4,1);
\node [below] at (3,1) {{\color{teal} $3A+B$}};
\node [left] at (2,2) {{\color{teal} $2A+2B$}};
\node [above] at (4,1) {{\color{teal} $4A+B$}};
\node [below] at (5,0) {{\color{teal} $5A$}};
\node [above] at (2.6,1.5) {{\color{teal} $k_3$}};
\node [below] at (4.5,0.5) {{\color{teal} $k_4$}};
\draw [->, line width=2, red] (3,1) -- (2,2);
\draw [->, line width=2, red] (4,1) -- (5,0);
\end{tikzpicture}
\caption{Reaction network embedded in Euclidean plane}
\end{subfigure}
\begin{subfigure}[b]{0.35\textwidth}
\includegraphics[scale=0.35]{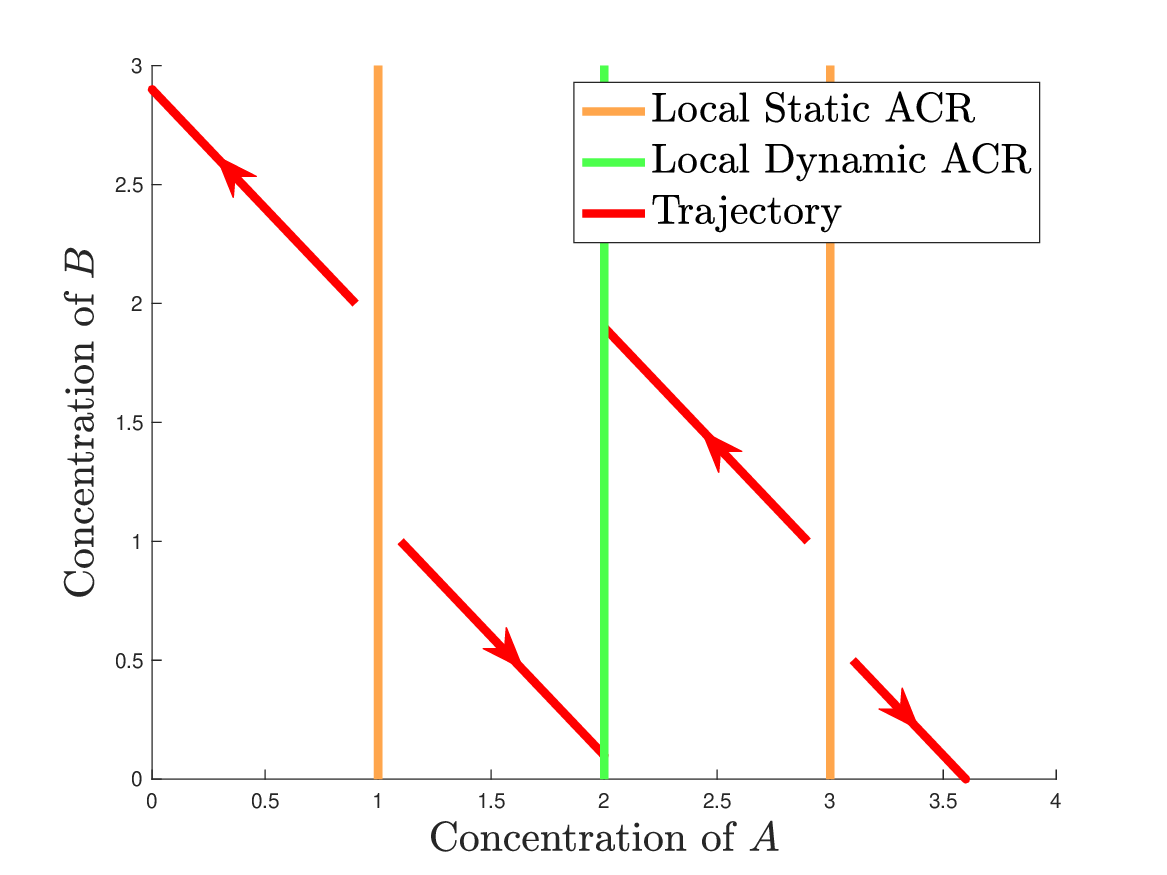}
\caption{Trajectories in phase plane}
\end{subfigure}
\caption{The mass action system of the reaction network ($A+B \xrightarrow{k_1} 2B,~ 2A+B \xrightarrow{k_2} 3A,~3A+B \xrightarrow{k_3} 2A+2B,~ 4A+B \xrightarrow{k_4} 5A$). For the rate constants $k_1=6, k_2=11,k_3=6, k_4=1$, $a$ is a local static and dynamic ACR variable.
There are three local static ACR values, $a=1, a=2$, and $a=3$. Only $a=2$ is a local dynamic ACR value.}
\label{fig:jeohijrtihj;rlh;jrp3op}
\end{figure}

\been
\item\label{ex:eogh3o5yih3} ({\it Almost cylinder ACR/Neighborhood ACR/Weak cylinder ACR}) Consider the reaction network shown in Figure \ref{fig:jeohijrtihj;rlh;jrp3op}(a). 
\been
\item For the mass action system defined by $k_1=6, k_2=11, k_3=6, k_4=1$, 
the positive steady states form three distinct rays $a=1$, $a=2$ and $a=3$. Therefore, $a$ is a cylinder static ACR variable with multiple ACR values $\{1,2,3\}$. It is easily checked (see simulated trajectories in Figure \ref{fig:jeohijrtihj;rlh;jrp3op}(b)), that only $a=2$ is locally stable within each compatibility class. The maximal cylinder which forms the basin of attraction for $\{a=2\}$ has radius $1$. 
\item For the mass action system defined by $k_1=1, k_2=3, k_3=3, k_4=1$, 
all positive steady states lie on $\{a=1\}$. Moreover, the positive steady states are repelling. It follows that $a$ is a (global) static ACR variable but not a local dynamic ACR variable. 
\item For the general case, the system of mass action ODEs is 
\alis{
\dot a = -ab(k_1 - k_2 a + k_3 a^2 - k_4 a^3), \quad \dot b = ab(k_1 - k_2 a + k_3 a^2 - k_4 a^3).
}
The univariate polynomial $k_1 - k_2 a + k_3 a^2 - k_4 a^3$ must have at least one positive, real zero. Moreover, when there are multiple positive zeros, say $a^*$ and $a^{**}$, clearly there are non-intersecting cylinders that contain the sets $\{a=a^*\}$ and $\{a=a^{**}\}$. This shows that $A$ is a local static ACR species, i.e. for any choice of mass action kinetics, $a$ is a local static ACR variable. In some cases, $a$ may be a global static ACR variable. However, $A$ is not a local dynamic ACR species as there is a choice of rate constants for which $a$ is not a local dynamic ACR variable. In other words, the reaction network has capacity for local dynamic ACR and is local static ACR  (see Definition \ref{def:049jhiwr} for meaning of `capacity'). 
\enen

\item ({\it Multiple ACR values in almost cylinder ACR}) Consider the reaction network shown in Figure \ref{fig:jeohijrtihj;rlh;jrp3op}(a). 
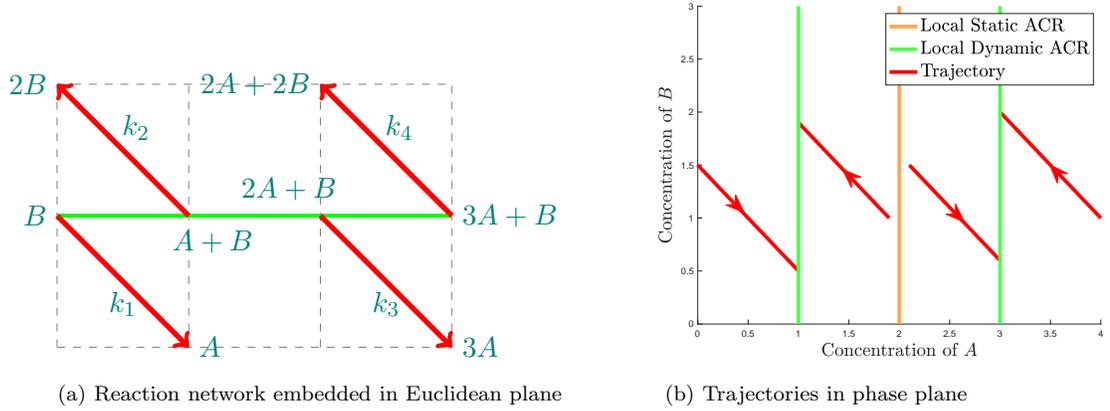
\begin{figure}[h!] 
\centering
\begin{subfigure}[b]{0.5\textwidth}
\begin{tikzpicture}[scale=1.75]
\draw[help lines, dashed, line width=0.25] (0,0) grid (3,2);

\node [right] at (1,0) {{\color{teal} $A$}};
\node [left] at (0,1) {{\color{teal} $B$}};
\node [right] at (0.8,0.8) {{\color{teal} $A+B$}};
\node [left] at (0,2) {{\color{teal} $2B$}};

\node [above] at (.6,1.5) {{\color{teal} $k_2$}};
\node [below] at (.5,0.5) {{\color{teal} $k_1$}};

\draw [-, line width=1.5, green] (0,1) -- (3,1);
\draw [->, line width=2, red] (0,1) -- (1,0);
\draw [->, line width=2, red] (1,1) -- (0,2);

\node [right] at (3,0) {{\color{teal} $3A$}};
\node [left] at (2.2,1.2) {{\color{teal} $2A+B$}};
\node [right] at (3,1) {{\color{teal} $3A+B$}};
\node [left] at (2,2) {{\color{teal} $2A+2B$}};

\node [above] at (2.6,1.5) {{\color{teal} $k_4$}};
\node [below] at (2.5,0.5) {{\color{teal} $k_3$}};

\draw [->, line width=2, red] (2,1) -- (3,0);
\draw [->, line width=2, red] (3,1) -- (2,2);

\end{tikzpicture}
\caption{Reaction network embedded in Euclidean plane}
\end{subfigure}
\begin{subfigure}[b]{0.3\textwidth}
\includegraphics[scale=0.35]{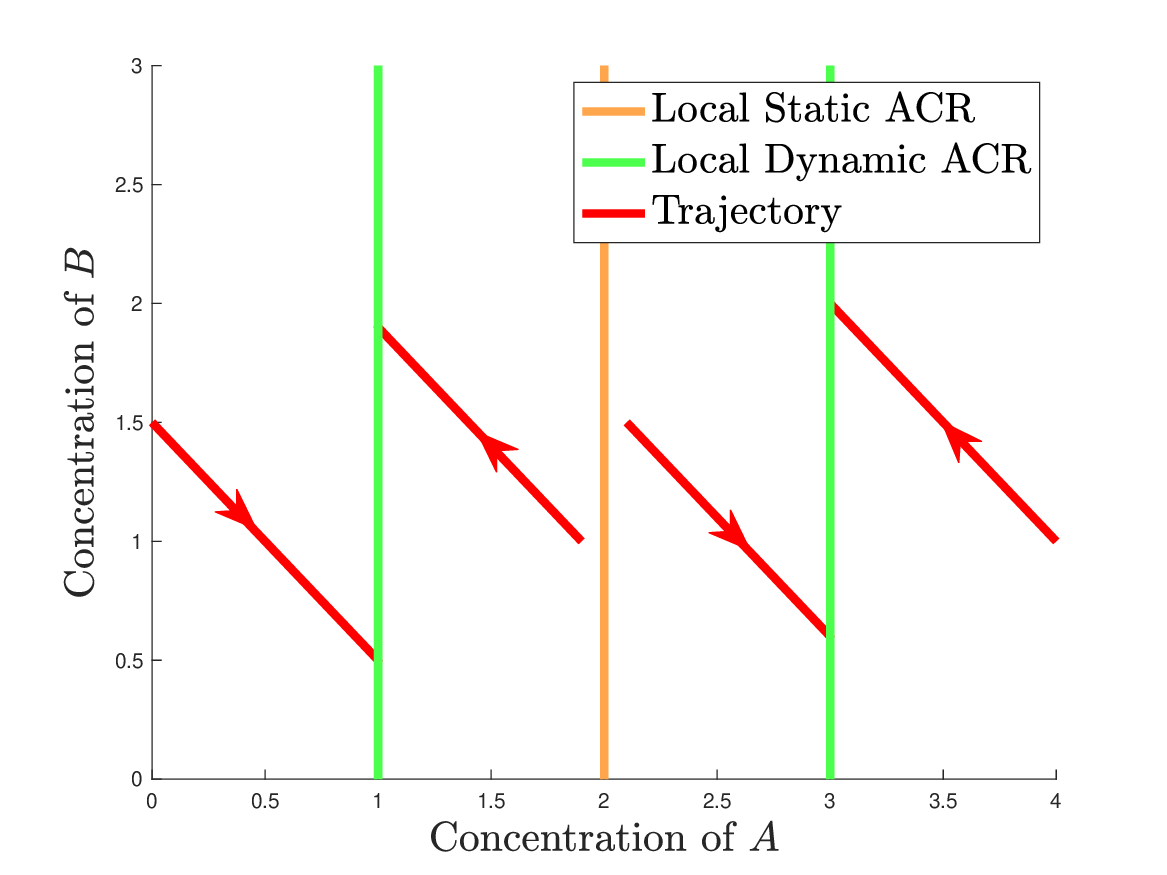}
\caption{Trajectories in phase plane}
\end{subfigure}
\caption{The mass action system of the reaction network ($A+B \xrightarrow{k_2} 2B,~ B \xrightarrow{k_1} A,~ 3A+B \xrightarrow{k_4} 2A+2B,~2A+B \xrightarrow{k_3} 3A$). For the rate constants $k_1=6, k_2=11,k_3=6, k_4=1$, $a$ is a local static and dynamic ACR variable.
There are three local static ACR values, $a=1, a=2$, and $a=3$, of which $a=1$ and $a=3$ are local dynamic ACR values.}
\label{fig:q4poiy5hoeghoiwfjoi2}
\end{figure}

%

\begin{figure}[h!] 
\centering
\includegraphics[scale=0.4]{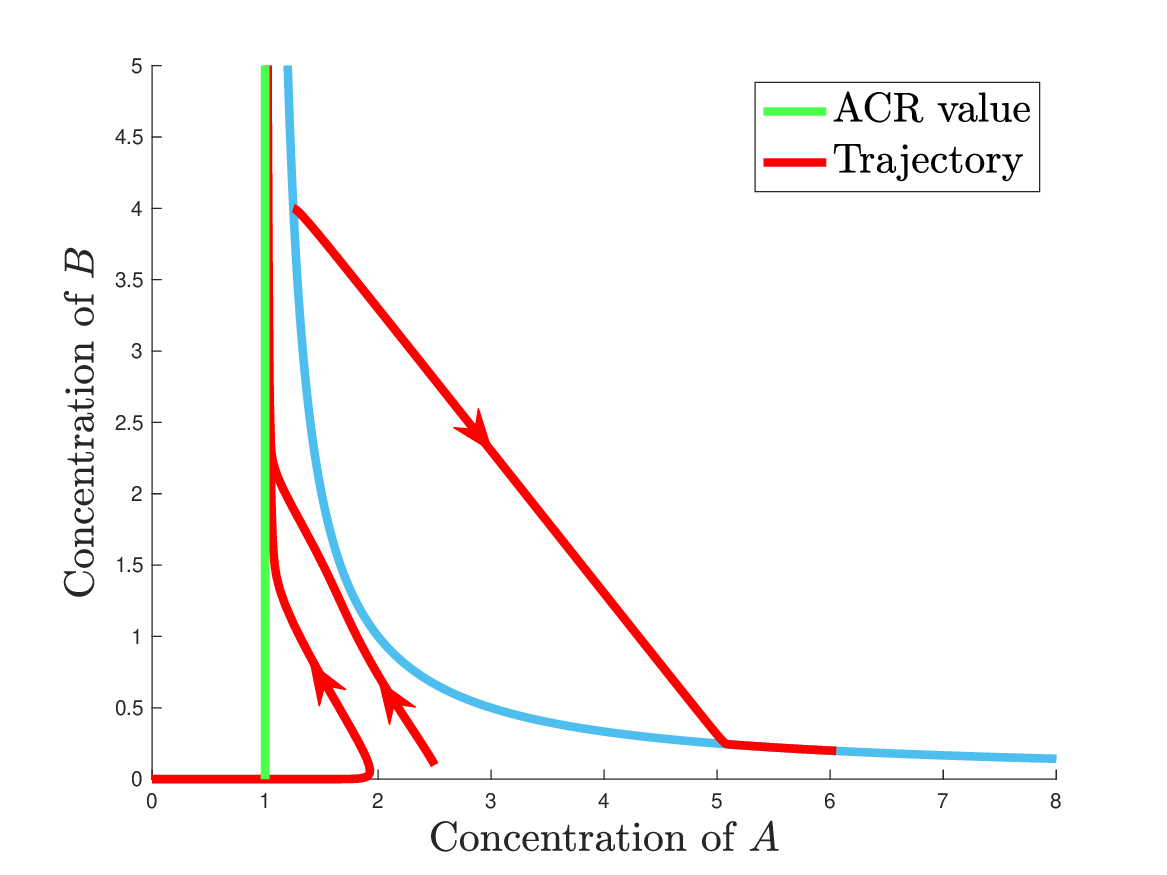}
\caption{(Neighborhood ACR but not Cylinder ACR:) Trajectories of the mass action system in \eqref{eq:qoh3oho35gwfjg}. An example of a mass action system which, in the concentration of $A$, (i) has neighborhood ACR as the strongest local ACR property, (ii) has half ACR as the strongest non-local ACR property, but (iii) does not have  cylinder ACR. For any initial value of $a(0) = 1+ \varepsilon$ for some $\varepsilon >0$, there is a $b(0)$ large enough such that the trajectory moves away from the ACR line. }
\label{fig:qoerho3gho3ghowg}
\end{figure}

\item ({\it Neighborhood ACR but not almost cylinder ACR}) 
Consider the mass action system shown below:
\begin{align} \label{eq:qoh3oho35gwfjg}
A+B \xrightarrow{1} 2A, &\quad 2A+B \stackrel[1]{1}{\rightleftarrows} A+2B, \nonumber\\
2A+2B \xrightarrow{2} A+3B, &\quad 3A+2B \xrightarrow{1} 4A+B.
\end{align}
The resulting mass action system when taken with inflows is 
\begin{align*}
\dot a &= ab (1-a)(1-(a-1)b) + g_a, \\
\dot b &= - ab (1-a)(1-(a-1)b) + g_b.
\end{align*}
Some trajectories are shown in Figure \ref{fig:qoerho3gho3ghowg} for the case of $g_a=0.1$ and $g_b=0$. Any initial value to the left of the ACR hyperplane converges to the ACR hyperplane, but for every initial $a$ value to the right of the ACR hyperplane, there is a large enough $b$ value such that the trajectory does not converge to the ACR value. 

\suspend{enumerate}

\subsection{An example of a system with subspace (but not full space) dynamic ACR}

The following example illustrates the need for defining subspace ACR. 

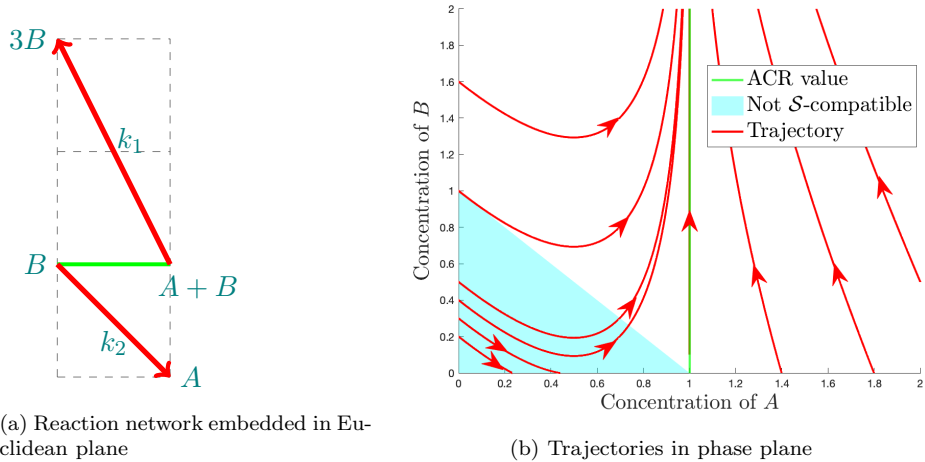
\begin{figure}[h!] 
\centering
\begin{subfigure}[b]{0.3\textwidth}
\begin{tikzpicture}[scale=1.5]
\draw[help lines, dashed, line width=0.25] (0,0) grid (1,3);

\node [right] at (1,0) {{\color{teal} $A$}};
\node [left] at (0,1) {{\color{teal} $B$}};
\node [right] at (0.8,0.8) {{\color{teal} $A+B$}};
\node [left] at (0,3) {{\color{teal} $3B$}};

\node [above] at (.65,1.9) {{\color{teal} $k_1$}};
\node [below] at (.5,0.5) {{\color{teal} $k_2$}};

\draw [-, line width=1.5, green] (0,1) -- (1,1);
\draw [->, line width=2, red] (0,1) -- (1,0);
\draw [->, line width=2, red] (1,1) -- (0,3);

\end{tikzpicture}
\caption{Reaction network embedded in Euclidean plane}
\end{subfigure}
\begin{subfigure}[b]{0.45\textwidth}
\includegraphics[scale=0.4]{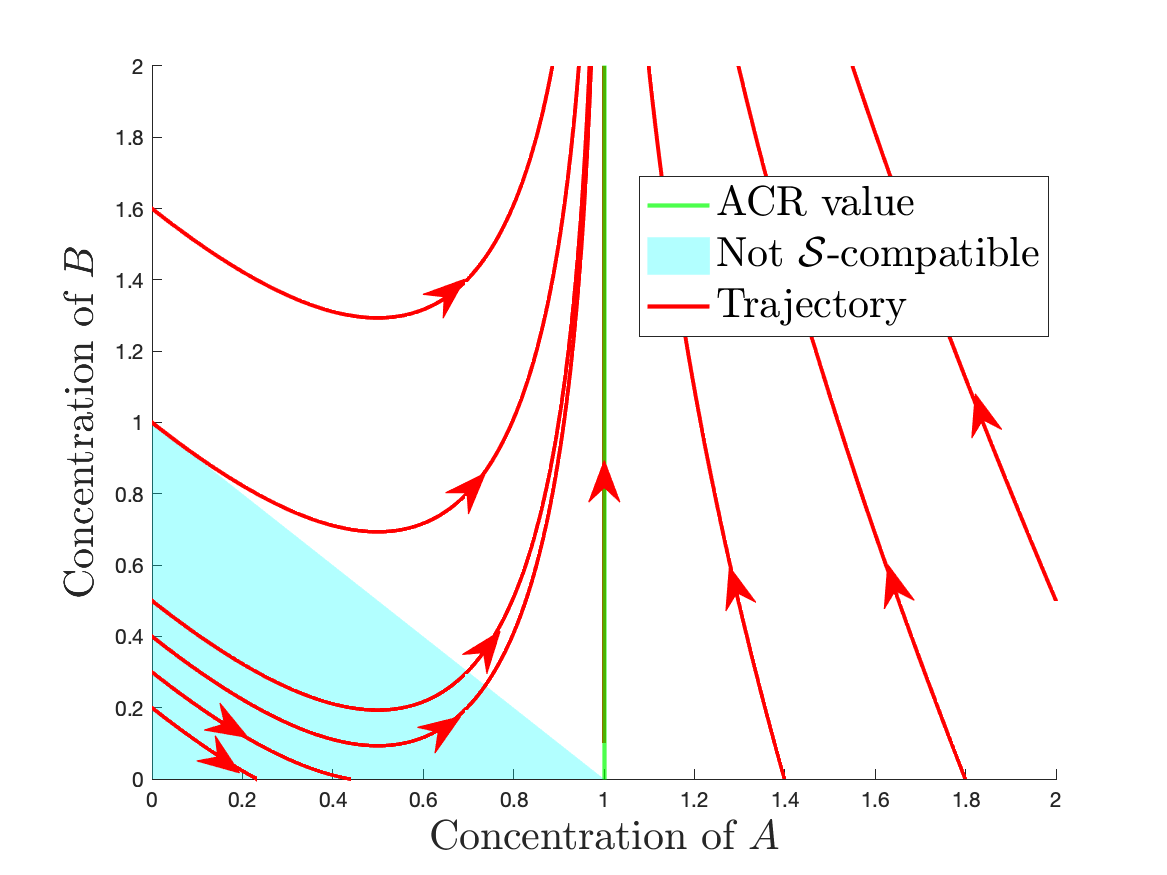}
\caption{Trajectories in phase plane}
\end{subfigure}
\caption{The mass action system $\{A+B \xrightarrow{k_1} 3B, B \xrightarrow{k_2} A\}$ is subspace dynamic ACR. For any choice of rate constants, all initial values that are $\SS$-compatible, where $\SS = \spn\{(-1,1)\}$ with the ACR hyperplane converge to the ACR hyperplane.}
\label{fig:;qerghoi3thoijfw}
\end{figure}

\resume{enumerate}
\item ({\it Subspace ACR and Cylinder  ACR}) Consider the following mass action system: 
\alis{
A+B \xrightarrow{k_1} 3B, \quad B \xrightarrow{k_2} A, 
}
which defines the ODE system:
\alis{
\dot a = b(k_2 -k_1 a), \quad \dot b = -b(k_2 - 2k_1 a). 
}
As long as $b(t)$ remains positive, $a(t)$ will move towards $a^* = k_2/k_1$. Since the stoichiometric subspace is all of $\R^2$, dynamic ACR requires every positive initial value to converge to $a^*$. This condition is not satisfied since, as trajectories in bottom left corner of Figure \ref{fig:;qerghoi3thoijfw}(b) show, there exist initial conditions which converge to the $b=0$ boundary and not to the ACR hyperplane. However, if we define $\Omega$ to be $\{a=a^*\} + \spn\{(-1,1)\}$, then $a$ is dynamic ACR w.r.t. $\Omega$. Thus, $a$ is subspace (dynamic) ACR. 

\enen

\section{Classification of Minimal Static and Dynamic ACR Networks} \label{sec:6yyj56hdgkwtk}

We consider networks of small size, ones with at most two reactions and at most two species. For such networks, we can catalogue many ACR properties. 
Some of these networks have archetypal ACR dynamics.
The study of minimal, archetypal motifs is valuable because it may reveal the underlying principles at play in the dynamics of more complex networks. 
\begin{definition} \label{def:049jhiwr}
Suppose that $(\GG,K)$ is a mass action system resulting from the reaction network $\GG$, where $K$ denotes the specific choice of mass action rate constants. 
Let $\PP \in \{$static, strong static, dynamic, weak dynamic, wide basin dynamic, narrow basin dynamic, full basin dynamic$\}$. 
\beit
\item We say that $\GG$ has {\em capacity for $\PP$-ACR} if there is a $K$ such that the mass action system $(\GG,K)$ is a $\PP$-ACR system.
\item We say that $\GG$ is a {\em $\PP$-ACR network} if $(\GG,K)$ is a $\PP$-ACR system for all choices of $K$. 
\item We say a species $X$ in a network $\GG$ is a {\em $\PP$-ACR species} if the concentration of $X$ is a $\PP$-ACR variable in $(\GG,K)$ for all choices of $K$. 
\enit
\end{definition}

\subsection{Static and dynamic ACR for reaction networks with only 1 reaction or only 1 species.}

\begin{theorem}[Static and Dynamic ACR in one-reaction networks]
A network with $n \ge 1$ species and only one reaction is neither static nor dynamic ACR for any choice of rate constants.  
\end{theorem}
\begin{proof}
A network with only one reaction has no positive steady states and is therefore not static ACR. 
Such a network is also not dynamic ACR since for each $i \in \{1, \ldots, n\}$, $\dot x_i$ is either strictly positive in the entire positive orthant, or strictly negative in the entire positive orthant, or identically zero in the entire positive orthant. Therefore, either $x_i$ goes to infinity, to zero, or $\dot x_i \equiv 0$. In every case $x_i$ fails to be a dynamic ACR variable. 
\end{proof}

\begin{theorem}[Static ACR in one-species networks]
Let $\GG$ be a network with $1$ species and arbitrary number of reactions. 
The following are equivalent.
\been[{A}1.]
\item $\GG$ has the capacity for static ACR. 
\item $(\GG,K)$ has a unique positive steady state for some $K$. 
\enen

The following are equivalent.
\been[{B}1.]
\item $\GG$ is static ACR. 
\item $(\GG,K)$ has a unique positive steady state for every choice of $K$. 
\enen
\end{theorem}
\begin{proof}
The results follow immediately from basic properties of one dimensional dynamical systems.
\end{proof}

\begin{theorem}[Dynamic ACR in one-species networks]
Let $\GG$ be a network with $1$ species and arbitrary number of reactions. 
The following are equivalent.
\been[{A}1.]
\item $\GG$ has the capacity for weak dynamic ACR. 
\item $\GG$ has the capacity for dynamic ACR. 
\item $(\GG,K)$ has a unique positive steady state for some $K$, and this steady state is globally attracting. 
\item $(\GG,K)$ is full basin dynamic ACR for some $K$. 
\enen

The following are equivalent.
\been[{B}1.]
\item $\GG$ is weak dynamic ACR. 
\item $\GG$ is dynamic ACR. 
\item For every choice of $K$, $(\GG,K)$ has a unique positive steady state, and this steady state is globally attracting. 
\item $(\GG,K)$ is full basin dynamic ACR for every $K$. 
\enen
\end{theorem}
\begin{proof}
The results follow immediately from basic properties of one dimensional dynamical systems.
\end{proof}
Many examples of one-species networks along with their ACR properties are shown in Table \ref{tab:5;yj4otihjo4jgjpe}. 

\begin{table}[htbp] 
\setlength\extrarowheight{0pt}
\centering
\begin{adjustbox}{max width=\textwidth}
 \begin{tabular}{|c|c||c|c||c|c|}
 \hline
 \specialcell{Network} & \specialcell{Arrow\\Diagram} & \specialcell{Capacity for\\static ACR?}& \specialcell{Is network\\static ACR?} & \specialcell{Capacity for\\dynamic ACR?}& \specialcell{Is network\\dynamic ACR?} \\
 \hline
 \hline
\specialcell{$0 \to A$} & $\longrightarrow$  & No & No  & No & No  \\
\hline
\specialcell{$0 \rlas A$} & $\longrightarrow, \longleftarrow$ & Yes & Yes  & Yes & Yes  \\
\hline
\specialcell{$0 \to A, 2A \to 3A$} & $\longrightarrow, \longrightarrow$ & No & No  & No & No   \\
\hline
\specialcell{$0 \to A, 2A \la 3A$} & $\longrightarrow, \longleftarrow$ & Yes & Yes  & Yes & Yes   \\
\hline
\specialcell{$0 \la A, 2A \to 3A$} & $\longleftarrow, \longrightarrow$ & Yes & Yes  & No & No   \\
\hline
\specialcell{$0 \la A, 2A \la 3A$} & $\longleftarrow, \longleftarrow$ & No & No  & No & No   \\
\hline
\specialcell{$2A \rlas 3A$} & $\longrightarrow, \longleftarrow$ & Yes & Yes  & Yes & Yes   \\
\hline
\specialcell{$0 \rlas A, 2A \to 3A$} & $\longrightarrow, \longleftarrow, \longrightarrow$ & Yes & No  & No & No   \\
\hline
\specialcell{$0 \rlas A, 2A \la 3A$} & $\longrightarrow, \longleftarrow, \longleftarrow$ & Yes & Yes  & Yes & Yes   \\
\hline
\specialcell{$0 \to A, 2A \rlas 3A$} & $\longrightarrow, \longrightarrow, \longleftarrow$ & Yes & Yes  & Yes & Yes   \\
\hline
\specialcell{$0 \la A, 2A \rlas 3A$} & $\longleftarrow, \longrightarrow, \longleftarrow$ & Yes & No  & No & No   \\
\hline
\specialcell{$0 \rlas A,~~ 2A \rlas 3A$} & $\longleftarrow, \longrightarrow, \longleftarrow, \longrightarrow$ & Yes & No  & Yes & No  \\
\hline
 \end{tabular}
 \end{adjustbox}
 \caption{Subnetworks of $0 \rlas A,~~ 2A \rlas 3A$ show diverse behaviors when catalogued according to capacities for static and dynamic ACR and according to whether the network is static or dynamic ACR. The diversity illustrates the range of possibilities even for one species networks.} \label{tab:5;yj4otihjo4jgjpe}
 \end{table}

 \begin{figure}[h!]
\centering
\begin{tikzpicture}[scale=1]
\draw[help lines, dashed, line width=0.25] (0,0) grid (6,6);

\node [below=0.25in] at (3,0) {{\color{olive} $X$}};
\node [left=0.25in] at (0,3) {{\color{olive} $Y$}};

\node [below left] at (2,3) {{\cbl $(a_1,b_1)$}};
\node [above right] at (5,2) {{\cbl $(a_2,b_2)$}};

\node [above right] at (4,5) {{\cbl $(\wt a_1,\wt b_1)$}};
\node [below left] at (1,1) {{\cbl $(\wt a_2,\wt b_2)$}};

\node [below] at (2,0) {{\color{teal} $a_1$}};
\node [below] at (4,0) {{\color{teal} $\wt a_1$}};

\node [left] at (0,3) {{\color{teal} $b_1$}};
\node [left] at (0,5) {{\color{teal} $\wt b_1$}};

\node [above left] at (3,4) {{\color{red} $k_1$}};
\node [below right] at (3,1.5) {{\color{red} $k_2$}};

\draw [->, line width=2, red] (2,3) -- (4,5);
\draw [->, line width=2, red] (5,2) -- (1,1);
\draw [-, line width=1.5, green] (2,3) -- (5,2);
\end{tikzpicture}
\vspace{-0.3cm}
\caption{A reaction network with two reactions and at most two species $X$ and $Y$ can be depicted as a pair of arrows embedded in the Euclidean plane $\R^2_{\ge 0}$. 
The red arrows depict reactions and the green line segment joining the two source complexes is the reactant polytope of a network with two reactions. 
The arrow from $(a_1,b_1)$ to $(\wt a_1,\wt b_1)$ portrays the reaction $a_1X + b_1 Y \to \wt a_1X + \wt b_1 Y$.  The label $k_1$ is the mass action rate constant of this reaction. 
The form of ACR, static or dynamic, as well as basin type can be decided based on the geometry of the three objects appearing in the Figure: the reactant polytope and the reaction arrows.}
\label{fig:4o268u460ujh4ij}
\end{figure}
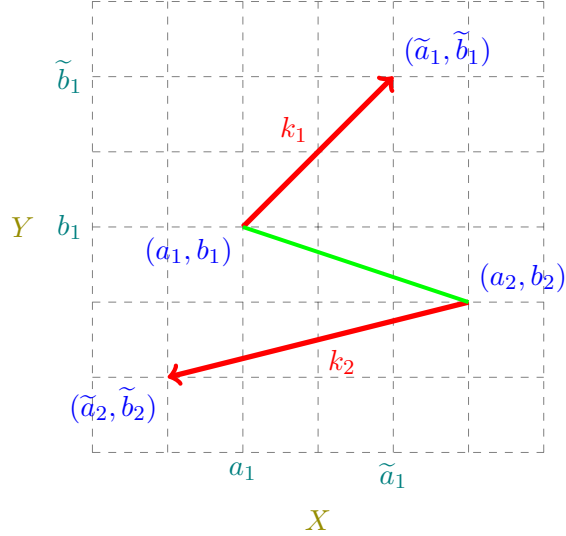

\subsection{Reaction networks with 2 reactions and 2 or fewer species: Notation}

We now classify reaction networks with two reactions and two species. We start by defining notation that will be used in the rest of the section. See Figure \ref{fig:4o268u460ujh4ij} for a geometric rendering of a reaction network, and how it relates to the notation. 
Let $\GG$ be a reaction network with at most 2 species ($X, Y$) and the following two reactions:
\begin{align} \label{eq:qpeogoe0thg}
a_1 X + b_1 Y \xrightarrow{k_1} \wt a_1 X + \wt b_1 Y, \quad  
a_2 X + b_2 Y \xrightarrow{k_2} \wt a_2 X + \wt b_2 Y, 
\end{align}
where $a_i, b_i, \wt a_i, \wt b_i \in \R_{\ge 0}$ and $(\wt a_i, \wt b_i) \ne (a_i, b_i)$ for  $i \in \{1,2\}$. 
Although stoichiometric coefficients are usually integers, we allow real values here since the results remain unchanged under this generality. 
The labels $k_1$ and $k_2$ are mass action reaction rate constants, and are therefore positive reals.
Let 
\begin{align} 
\SS = \spn\left\{v_1 \coloneqq \begin{pmatrix}\wt a_1-a_1\\\wt b_1-b_1\end{pmatrix}, v_2 \coloneqq \begin{pmatrix}\wt a_2-a_2\\ \wt b_2-b_2\end{pmatrix}\right\}
\end{align}
 be the stoichiometric subspace of $\GG$. The mass action dynamical system $(\GG,K)$ explicitly is: 
\begin{align} \label{eq:34oyhhfiwjfp2j}
\dot x &= k_1(\wt a_1 - a_1) x^{a_1} y^{b_1} + k_2(\wt a_2 - a_2) x^{a_2} y^{b_2} \nonumber \\ 
\dot y &= k_1(\wt b_1 - b_1) x^{a_1} y^{b_1} + k_2(\wt b_2 - b_2) x^{a_2} y^{b_2}. 
\end{align}

\subsection{Static ACR for reaction networks with 2 reactions and 2 or fewer species.}
\label{sec:staticACR}

\begin{theorem}[Static ACR in networks with 2 reactions and 2 or fewer species] \label{thm:po43ghi3lwkfh}
Let $\GG$ be as in \eqref{eq:qpeogoe0thg}-\eqref{eq:34oyhhfiwjfp2j}. The following are equivalent: 
\been
\item $\GG$ has the capacity for static ACR. 
\item $\GG$ is static ACR. 
\item $\GG$ is strong static ACR.
\item 
\beit 
\item the two source complexes are different:
$
(a_1,b_1) \ne (a_2,b_2), 
$
\item the source complexes share a common coordinate:
$(a_2-a_1)(b_2 - b_1)=0$, and 
\item reaction vectors are negative scalar multiples of each other:
$
v_1 = - \mu v_2 \mbox{ for some } \mu > 0, 
$
(in particular $\dim(\SS)=1$). 
\enit
\enen
Furthermore, the following hold when $\GG$ has 2 species and is static ACR.
\beit
\item Either $X$ or $Y$, but not both, is a static ACR species. 
\item $X$ is an ACR species if $a_2 \ne a_1$. The variable $x$ has the static ACR value $\ds \left(k_2/(\mu k_1) \right)^{\frac{1}{a_1-a_2}}$. 
\item $Y$ is an ACR species if $b_2 \ne b_1$. The variable $y$ has the static ACR value $\ds \left(k_2/(\mu k_1) \right)^{\frac{1}{b_1-b_2}}$. 
\enit

\end{theorem}

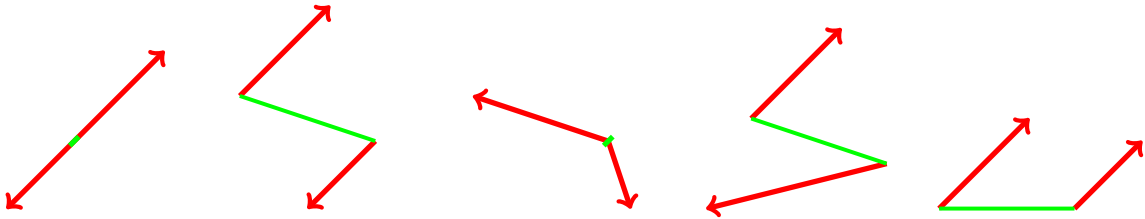
\begin{figure}[h!] 
\centering
\begin{subfigure}[b]{0.18\textwidth}
\begin{tikzpicture}[scale=0.6]
\draw [->, line width=2, red] (2,3) -- (4,5);
\draw [->, line width=2, red] (2,3) -- (0.5,1.5);
\draw [-, line width=2, green] (1.9,2.9) -- (2.1,3.1);
\end{tikzpicture}
\end{subfigure}
\begin{subfigure}[b]{0.18\textwidth}
\begin{tikzpicture}[scale=0.6]
\draw [->, line width=2, red] (2,3) -- (4,5);
\draw [->, line width=2, red] (5,2) -- (3.5,0.5);
\draw [-, line width=1.5, green] (2,3) -- (5,2);
\end{tikzpicture}
\end{subfigure}
\begin{subfigure}[b]{0.18\textwidth}
\begin{tikzpicture}[scale=0.6]
\draw [->, line width=2, red] (2,3) -- (-1,4);
\draw [->, line width=2, red] (2,3) -- (2.5,1.5);
\draw [-, line width=2, green] (1.9,2.9) -- (2.1,3.1);
\end{tikzpicture}
\end{subfigure}
\begin{subfigure}[b]{0.18\textwidth}
\begin{tikzpicture}[scale=0.6]
\draw [->, line width=2, red] (2,3) -- (4,5);
\draw [->, line width=2, red] (5,2) -- (1,1);
\draw [-, line width=1.5, green] (2,3) -- (5,2);
\end{tikzpicture}
\end{subfigure}
\begin{subfigure}[b]{0.18\textwidth}
\begin{tikzpicture}[scale=0.6]
\draw [->, line width=2, red] (2,3) -- (4,5);
\draw [->, line width=2, red] (5,3) -- (6.5,4.5);
\draw [-, line width=1.5, green] (2,3) -- (5,3);
\end{tikzpicture}
\end{subfigure}
\caption{Motifs that do not have the capacity for static ACR. A reaction network does not have the capacity for static ACR if the source complexes of the two reactions are the same, or if both  coordinates of the source complexes are different, or if the reaction vectors do not point in opposite directions. See Theorem \ref{thm:po43ghi3lwkfh} for precise conditions.}
\label{fig:qo53yihqengk3r}
\end{figure}

\begin{proof}
(${\it 3} \implies {\it 2} \implies {\it 1}$) holds by definition. 
We now show that (${\it 1} \implies {\it 4}$). 
Suppose that $v_1 \ne - \mu v_2$ for any $\mu > 0$, then there are no positive steady states for any choice of mass action rate constants and so $\GG$ does not have the capacity for static ACR. 

Now assume that $v_1 = -\mu v_2$ for some $\mu >0$. Then the mass action system is
\begin{align} \label{eq:owetihowh}
\dot x = (\wt a_1 - a_1) \left(k_1 x^{a_1}y^{b_1} - \frac1{\mu}k_2 x^{a_2} y^{b_2}\right), \quad \dot y = (\wt b_1 - b_1) \left(k_1 x^{a_1}y^{b_1} - \frac1{\mu}k_2 x^{a_2} y^{b_2}\right). 
\end{align}
If $(a_1,b_1) = (a_2,b_2)$, then 
\begin{align*}
\dot x = (\wt a_1 - a_1) \left(k_1  - \frac1{\mu}k_2 \right) x^{a_1}y^{b_1}, \quad 
\dot y = (\wt b_1 - b_1) \left(k_1  - \frac1{\mu}k_2 \right) x^{a_1}y^{b_1}.
\end{align*} 
If $k_1=k_2/\mu$, then every positive point is a steady state and so the system is not static ACR. If $k_1 \ne k_2/\mu$, then at least one of $\dot x$ or $\dot y$ is either positive on all of $\R^2_{>0}$ or negative on all of $\R^2_{>0}$. But then there is no positive steady state. So $\GG$ does not have the capacity for static ACR. 

From \eqref{eq:owetihowh}, steady states must satisfy the equation
\begin{align} \label{eq:3oth3tugh}
x^{a_1-a_2} y^{b_1-b_2} = \frac{k_2}{\mu k_1} =: k. 
\end{align}
Now suppose that $0 \not \in \{a_2-a_1, b_2 - b_1\}$. From \eqref{eq:3oth3tugh}, we see that 
\[
(x_\beta,y_\beta) = \left( (k\beta)^{1/(a_1-a_2)}, \beta^{-1/(b_1-b_2)} \right)
\]
is a steady state for every $\beta \in \R_{> 0}$. In particular, two distinct choices of $\beta$ result in distinct $x$ and $y$ components in the two steady states. This implies that neither variable is static ACR for any $k$. So $\GG$ does not have the capacity for static ACR. 

Finally, we show that (${\it 4} \implies {\it 3}$). 
Assume without loss of generality that $a_2 - a_1 \ne 0$ and $b_2 - b_1=0$. From \eqref{eq:3oth3tugh}, we have that 
$
x = k^* \coloneqq k^{\frac{1}{a_1-a_2}}
$ at steady state, which shows that the system is static ACR and $x$ is the static ACR variable. Since this is true for every choice of mass action rate constants, $\GG$ is static ACR and $X$ is a static ACR species. Every point on the hyperplane $\HH[x,k^*]$ is a steady state which shows that $X$ is strong static ACR. 

It is clear that the roles of species $X$ and $Y$ are reversed if we assume that $a_2=a_1$ and $b_2 \ne b_1$, which proves the claims about the species $Y$. 
\end{proof}

Conditions for static ACR in reaction networks with 2 reactions and 2 or fewer species have also been studied in \cite{meshkat2021absolute}.




\subsection{Network motifs and their embeddings} \label{sec:networkmotifs}

\begin{figure}[h!] 
\centering
\begin{subfigure}[b]{0.1\textwidth}
\begin{tikzpicture}[scale=1]
\draw [-, line width=2, green] ({0-1},{0}) -- ({0+1},{0});
\draw [->, line width=2, red] ({0-1},{0}) -- ({0-0.25},{0+1.5});
\draw [->, line width=2, red] ({0+1},{0}) -- ({0+0.25},{0+1.5});
\end{tikzpicture}
\caption{Motif}
\end{subfigure}

\begin{subfigure}[b]{0.25\textwidth}
\begin{tikzpicture}[scale=1.5]
\draw[help lines, dashed, line width=0.25] (0,0) grid (1,2);

\node [right] at (1,0) {{\color{teal} $A$}};
\node [left] at (0,2) {{\color{teal} $2B$}};
\node [right] at (1,2) {{\color{teal} $A+2B$}};
\node [left] at (0,0) {{\color{teal} $0$}};

\draw [-, line width=1.5, green] (0,0) -- (1,0);
\draw [->, line width=2, red] (0,0) -- (1,2);
\draw [->, line width=2, red] (1,0) -- (0,2);

\end{tikzpicture}
\caption{Embedding 1}
\end{subfigure}
\quad \quad \quad \quad
\begin{subfigure}[b]{0.25\textwidth}
\begin{tikzpicture}[scale=1.5]
\draw[help lines, dashed, line width=0.25] (0,0) grid (3,2);

\node [right] at (3,0) {{\color{teal} $3A$}};
\node [right] at (1.1,2) {{\color{teal} $A+2B$}};
\node [left] at (.9,1) {{\color{teal} $A+B$}};
\node [left] at (0,0) {{\color{teal} $0$}};

\draw [-, line width=1.5, green] (0,0) -- (3,0);
\draw [->, line width=2, red] (0,0) -- (1,1);
\draw [->, line width=2, red] (3,0) -- (1,2);


\end{tikzpicture}
\caption{Embedding 2}
\end{subfigure}
\caption{A network motif (on top) and two of its embeddings (bottom row).}
\label{fig:o3gh3oigjwwpg}
\end{figure}
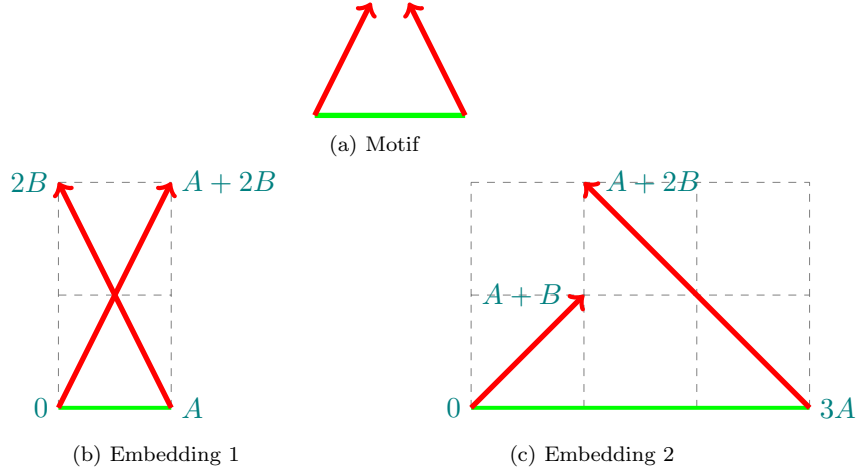

Similar to static ACR, we will show that dynamic ACR is a network property. If a network has the capacity for dynamic ACR, then it is dynamic ACR. Moreover, whether a network has the capacity for dynamic ACR depends only on its topology and not on the specific embedding in the Euclidean plane. We refer to such a class  of networks as a motif. 
A network motif in two dimensions is determined by: (i) slope of the reactant polytope, (ii) the quadrant or axis each reaction points along, and (iii) the relative slopes of the two reactions. 
We demonstrate the relation between a network motif and its multiple embeddings via an example in Figure \ref{fig:o3gh3oigjwwpg}. 

\subsection{Dynamic ACR in static ACR networks with 2 reactions \& 2 or fewer species}

\begin{theorem}[Dynamic ACR in static ACR networks with 2 reactions \& 2 or fewer species] \label{thm:oh3ihgj3oth3}
Let $\GG$ be as in \eqref{eq:qpeogoe0thg}-\eqref{eq:34oyhhfiwjfp2j} and suppose that $\dim(\SS)=1$.  The following statements {\it A}1-{\it A}6 are equivalent: 
\been[{A}1.]
\item $\GG$ has the capacity for weak dynamic ACR. 
\item $\GG$ is weak dynamic ACR. 
\item $\GG$ is full basin weak dynamic ACR. 
\item $\GG$ has the capacity for dynamic ACR. 
\item $\GG$ is (full space) dynamic ACR. 
\item $\GG$ is static ACR and $(\wt a_1-a_1)(a_2-a_1) + (\wt b_1-b_1)(b_2-b_1) > 0$. 
\enen
Now, suppose that $\GG$ has 2 species and is dynamic ACR. Then either $X$ or $Y$, but not both, is a dynamic ACR species. Furthermore, the following statements {\it B}1-{\it B}3 are equivalent. 
\been[{B}1.]
\item $X$ is a dynamic ACR species. 
\item $X$ is a static ACR species. 
\item $a_2 \ne a_1$. 
\enen
\beit
\item When $X$ is a dynamic ACR species, \\
\[
\mbox{the dynamic ACR value of the variable } x = \mbox{ the static ACR value of } x = \left( \frac{k_2}{\mu k_1} \right)^{\frac{1}{a_1-a_2}}. 
\]
\enit
Analogous statements to {\it B}1-{\it B}3 and the statement about ACR value hold when $X$ is replaced with $Y$, and $a_i$ is replaced with $b_i$ for $i \in \{1,2\}$. 

\end{theorem}
\begin{proof}
(${\it A3} \implies {\it A2} \implies {\it A1}$) and (${\it A5} \implies {\it A4}\implies {\it A1}$) hold by definition.  

We now show that (${\it A1} \implies {\it A6} \implies {\it A5, A3}$). Suppose that $\GG$ has the capacity for weak dynamic ACR. From properties of one dimensional dynamical systems, it is clear that $\GG$ has the capacity for static ACR. From Theorem \ref{thm:po43ghi3lwkfh}, $\GG$ is static ACR and one of the two (but not both) is a static ACR species.
Assume that $X$ (and not $Y$) is the static ACR species. Then $b_2=b_1$, $a_2 \ne a_1$, and $x$ is static ACR with value $\ds x^* \coloneqq \left( k_2/(\mu k_1) \right)^{\frac{1}{a_1-a_2}}$. 
\begin{align*}
\dot x = (\wt a_1 - a_1) y^{b_1} x^{a_1} \left(k_1  - \frac1{\mu}k_2 x^{a_2 - a_1} \right). 
\end{align*}
Clearly, the steady state $x^*$ is stable if and only if $(\wt a_1 - a_1)(a_2-a_1) > 0$. If we assume instead that $Y$ (and not $X$) is the static ACR species, we get the stability condition $(\wt b_1 - b_1)(b_2-b_1) > 0$. The desired inequality in {\it A6} is obtained by combining the two stability conditions, since it is always the case that one term in {\it A6} is positive and the other term is zero which shows that (${\it A1} \implies {\it A6}$). Since a unique (within a compatibility class) steady state that is stable must be globally stable for a one dimensional system (i.e. attracts all compatible positive points),  we also have that (${\it A6} \implies {\it A5}$). Moreover, the initial values that are not compatible with the hyperplane of steady states $\{x=x^*\}$ also result in trajectories that  move towards $\{x=x^*\}$ but converge at a boundary steady state. This gives full basin weak dynamic ACR, so we have also proved that  (${\it A6} \implies {\it A3}$).

The last part also shows that when $\GG$ is dynamic ACR, ({\it B3} $\implies$ {\it B1} $\implies$ {\it B2}) and ({\it B2} $\implies$ {\it B3}) is from Theorem \ref{thm:po43ghi3lwkfh}. The statement about the dynamic ACR value of $x$ also follows from the last part. Clearly, the assumption that $a_2=a_1$ and $b_2 \ne b_1$ will switch the roles of $X$ and $Y$ and give analogous statement for species $Y$.
\end{proof}

\begin{figure}[h!] 
\centering
\begin{tikzpicture}[scale=1]
\draw [->, line width=2, red] ({4*cos(90)+0.5},{4*sin(90)}) -- ({4*cos(90)+0.5},{4*sin(90)-1});
\draw [->, line width=2, red] ({4*cos(90)-0.5},{4*sin(90)}) -- ({4*cos(90)-0.5},{4*sin(90)+1});
\draw [-, line width=2, green] ({4*cos(90)-0.5},{4*sin(90)}) -- ({4*cos(90)+0.5},{4*sin(90)});

\draw [->, line width=2, red] ({4*cos(45)+0.5},{4*sin(45)}) -- ({4*cos(45)-0.5},{4*sin(45)-1});
\draw [->, line width=2, red] ({4*cos(45)-0.5},{4*sin(45)}) -- ({4*cos(45)+0.5},{4*sin(45)+1});
\draw [-, line width=2, green] ({4*cos(45)-0.5},{4*sin(45)}) -- ({4*cos(45)+0.5},{4*sin(45)});

\draw [->, line width=2, red] ({4*cos(0)+0.5},{4*sin(0)}) -- ({4*cos(0)-0.5},{4*sin(0)});
\draw [->, line width=2, red] ({4*cos(0)-0.5},{4*sin(0)}) -- ({4*cos(0)+0.5},{4*sin(0)});
\draw [-, line width=2, green] ({4*cos(0)-0.5},{4*sin(0)}) -- ({4*cos(0)+0.5},{4*sin(0)});

\draw [->, line width=2, red] ({4*cos(-45)+0.5},{4*sin(-45)}) -- ({4*cos(-45)-0.5},{4*sin(-45)+1});
\draw [->, line width=2, red] ({4*cos(-45)-0.5},{4*sin(-45)}) -- ({4*cos(-45)+0.5},{4*sin(-45)-1});
\draw [-, line width=2, green] ({4*cos(-45)-0.5},{4*sin(-45)}) -- ({4*cos(-45)+0.5},{4*sin(-45)});

\draw [->, line width=2, red] ({4*cos(-90)+0.5},{4*sin(-90)}) -- ({4*cos(-90)+0.5},{4*sin(-90)+1});
\draw [->, line width=2, red] ({4*cos(-90)-0.5},{4*sin(-90)}) -- ({4*cos(-90)-0.5},{4*sin(-90)-1});
\draw [-, line width=2, green] ({4*cos(-90)-0.5},{4*sin(-90)}) -- ({4*cos(-90)+0.5},{4*sin(-90)});

\draw [->, line width=2, red] ({4*cos(-135)+0.5},{4*sin(-135)}) -- ({4*cos(-135)+1.5},{4*sin(-135)+1});
\draw [->, line width=2, red] ({4*cos(-135)-0.5},{4*sin(-135)}) -- ({4*cos(-135)-1.5},{4*sin(-135)-1});
\draw [-, line width=2, green] ({4*cos(-135)-0.5},{4*sin(-135)}) -- ({4*cos(-135)+0.5},{4*sin(-135)});

\draw [->, line width=2, red] ({4*cos(-135)+0.5},{4*sin(-135)}) -- ({4*cos(-135)+1.5},{4*sin(-135)+1});
\draw [->, line width=2, red] ({4*cos(-135)-0.5},{4*sin(-135)}) -- ({4*cos(-135)-1.5},{4*sin(-135)-1});
\draw [-, line width=2, green] ({4*cos(-135)-0.5},{4*sin(-135)}) -- ({4*cos(-135)+0.5},{4*sin(-135)});

\draw [->, line width=2, red] ({4*cos(180)+0.5},{4*sin(180)}) -- ({4*cos(180)+1.5},{4*sin(180)});
\draw [->, line width=2, red] ({4*cos(180)-0.5},{4*sin(180)}) -- ({4*cos(180)-1.5},{4*sin(180)});
\draw [-, line width=2, green] ({4*cos(180)-0.5},{4*sin(180)}) -- ({4*cos(180)+0.5},{4*sin(180)});

\draw [->, line width=2, red] ({4*cos(135)+0.5},{4*sin(135)}) -- ({4*cos(135)+1.5},{4*sin(135)-1});
\draw [->, line width=2, red] ({4*cos(135)-0.5},{4*sin(135)}) -- ({4*cos(135)-1.5},{4*sin(135)+1});
\draw [-, line width=2, green] ({4*cos(135)-0.5},{4*sin(135)}) -- ({4*cos(135)+0.5},{4*sin(135)});

\draw[line width=1, teal] ({3.5*cos(0)},{3.5*sin(0)}) ellipse (1.5cm and 1cm);
\draw[line width=1, cyan] ({3.75*cos(-22.5)},{3.75*sin(-22.5)}) ellipse (2cm and 3cm);
\draw[line width=1, magenta] ({3.5*cos(0)},{3.5*sin(0)}) ellipse (2.5cm and 5cm);
\draw[line width=1, olive] (0,0) ellipse (7.5cm and 6cm);

\node[shape=ellipse, draw=none, line width=0]() at (4.5,2.5) {\color{magenta} \specialcell{Dynamic\\ACR}};
\node[shape=ellipse, draw=none, line width=0]() at (4.25,-2.5) {\color{cyan} \specialcell{Wide basin\\dynamic\\ACR}};
\node[shape=ellipse, draw=none, line width=0]() at (3.5,-0.5) {\color{teal} \specialcell{Full basin}};
\node[shape=ellipse, draw=none, line width=0]() at (-0.75,0) {\color{olive} \specialcell{\Large Static\\\Large ACR\\\Large Motifs}};
\end{tikzpicture}
\caption{Network motifs with 2 reactions and 1 dimensional stoichiometric subspace that show static ACR. The three network types in the {\cma magenta} ellipse are dynamic ACR, the two network types in the smaller {\ccy cyan} ellipse are wide basin dynamic ACR, while the one network type in the smallest {\cte teal} ellipse is full basin dynamic ACR. In each of the networks, the static ACR species is $X$, which is on the axis parallel to the reaction polytope.}
\label{fig:o45yuohdaguh3iur}
\end{figure}
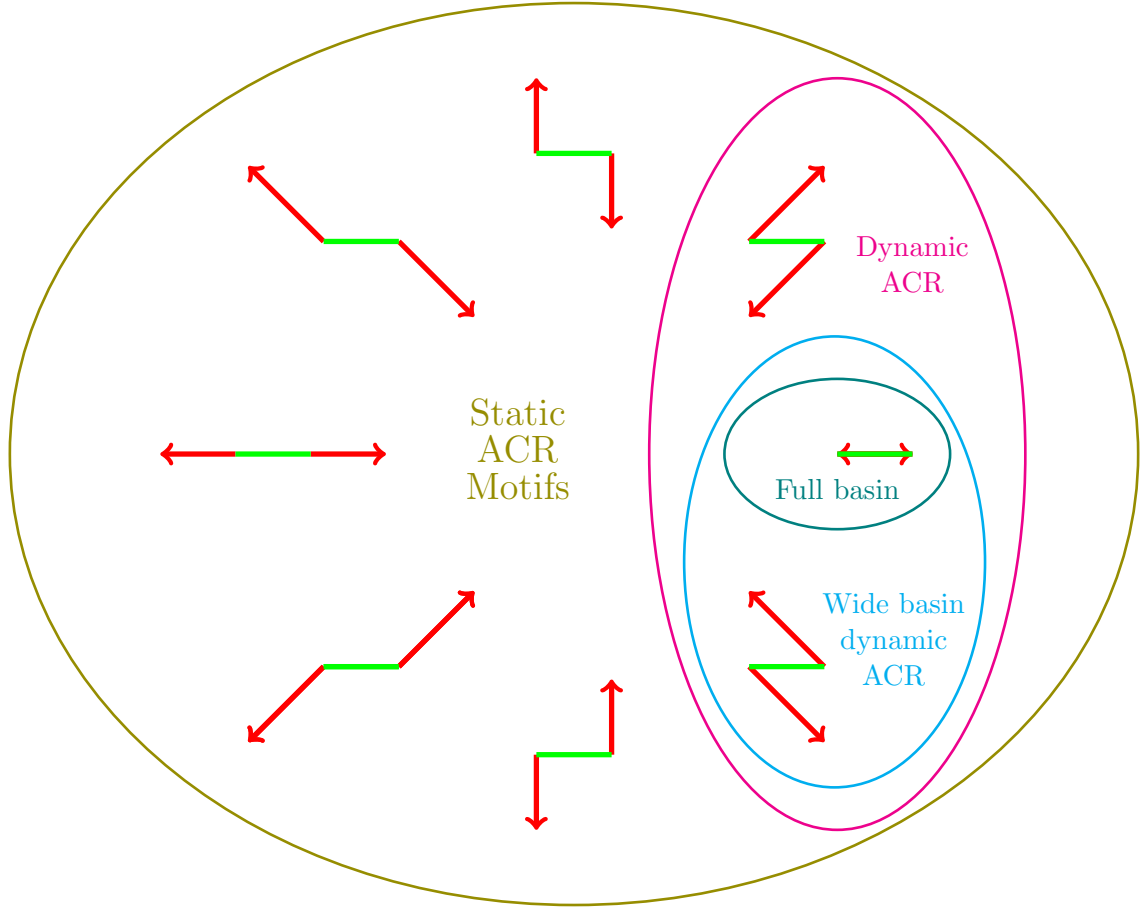

\begin{theorem}[Wide basin dynamic ACR in 2 reaction, 2 or fewer species networks with $\dim(\SS)=1$] \label{thm:31oy13thqog}
Let $\GG$ be as in \eqref{eq:qpeogoe0thg}-\eqref{eq:34oyhhfiwjfp2j}. Suppose that $\dim(\SS)=1$ and $\GG$ is dynamic ACR. 
\been
\item The dynamic ACR species is wide basin if and only if $(\wt a_1-a_1)(\wt b_1-b_1) \le 0$. 
\item The dynamic ACR species is full basin if and only if $(\wt a_1-a_1)(\wt b_1-b_1) = 0$. 
\enen
\end{theorem}
\begin{proof}
Since $\GG$ is dynamic ACR and $\dim(\SS)=1$, by Theorem \ref{thm:oh3ihgj3oth3}, $\GG$ is static ACR. By Theorem \ref{thm:po43ghi3lwkfh}, $v_1 = -\mu v_2$ for some $\mu >0$, so that the mass action system can be written as in \eqref{eq:owetihowh}.
This implies that $(\wt b_1 - b_1) \dot x - (\wt a_1 - a_1) \dot y=0$, and so two points $(x^*, y^*)$ and $(x_0, y_0)$ are compatible if and only if 
\begin{align}\label{eq:eoihnoigho3}
&(\wt b_1 - b_1) x^* - (\wt a_1 - a_1) y^*=(\wt b_1 - b_1) x_0 - (\wt a_1 - a_1) y_0 \nonumber \\ 
&\iff (\wt b_1 - b_1) (x_0 - x^*) =  (\wt a_1 - a_1) (y_0 - y^*)
\end{align}
Since at least one of $\wt a_1 - a_1$ or $\wt b_1 - b_1$ must be nonzero, the above implies that $\sgn((x_0 - x^*)(y_0 - y^*)) = \sgn((\wt b_1 - b_1)(\wt a_1 - a_1))$, where $\sgn$ is the sign function that has range $\{+, -, 0\}$ and is defined via:
\begin{align*}
\sgn(z) = \begin{cases}
+ &\mbox{ if } z > 0\\
- &\mbox{ if } z < 0\\
0 &\mbox{ if } z = 0
\end{cases}
\end{align*} 
We will assume without loss of generality that $X$ is the dynamic ACR species, so that for a particular choice of $K$, $x$ is the dynamic ACR variable with ACR value $k^*$ that depends on $K$. 

Suppose that $\sgn((\wt b_1 - b_1)(\wt a_1 - a_1)) = +$. Then $\alpha \coloneqq (\wt a_1 - a_1)/(\wt b_1 - b_1) > 0$. From \eqref{eq:eoihnoigho3}, $(x_0,y_0)$ is compatible with $(x^*,y^*)$ if and only if $x_0 - x^* = \alpha(y_0 - y^*)$ with $\alpha > 0$. 
Any initial value $(x_0,y_0)$ with $x_0 > k^* + \alpha y_0$ is incompatible with $\{x=k^*\}$.
To see this, note that 
\[
x_0 - k^* > \alpha y_0> \alpha(y_0 - y^*).
\]
Clearly $x_0$ is not bounded above on the incompatible set $\{x_0 > k^* + \alpha y_0\}$, so in this case $x$ is a narrow basin dynamic ACR variable for any choice of $K$. In other words, $X$ is a narrow basin dynamic ACR species. 

Now suppose that $\sgn((\wt b_1 - b_1)(\wt a_1 - a_1)) = 0$. Since $x$ is the dynamic ACR variable, $b_2 = b_1$. By condition in the previous part, $\wt a_1 \ne a_1$ and so $\wt b_1 = b_1$. In particular $\dot y =0$ and every positive initial value is compatible with the set $\{x = k^*\}$. So $x$ is a full basin dynamic ACR variable for every $K$. In other words, $X$ is a full basin dynamic ACR species.

Finally, suppose that $\sgn((\wt b_1 - b_1)(\wt a_1 - a_1)) = -$. In this case, $x_0 - x^* = - \beta (y_0 - y^*)$ with $\beta = -(\wt a_1 - a_1)/(\wt b_1 - b_1) > 0$. Any initial value $(x_0,y_0)$ with $x_0 + \beta y_0 > k^*$ is compatible with $\{x=k^*\}$. To see this, let 
\[
y^* \coloneqq (x_0 + \beta y_0 - k^*)/\beta >0.
\]
Then $(x_0,y_0)$ is compatible with $(k^*,y^*)$. So $x$ is a wide basin dynamic ACR variable for every $K$. In other words, $X$ is a wide basin dynamic ACR species. 
The only remaining thing to show is that $X$ is not full basin. Let $(x_0,y_0) \coloneqq \left(k^*/4, k^*/(4\beta) \right) \in \R^2_{>0}$. Then 
\[
x_0 + \beta y_0 = k^*/2 < k^* < k^* + \beta y^*, 
\]
and so $(x_0,y_0)$ is not compatible with $\{x = k^*\}$. 

This completes the proof since we have covered the entire range of $\sgn((\wt b_1 - b_1)(\wt a_1 - a_1))$. 
\end{proof}

The results of Theorems \ref{thm:po43ghi3lwkfh}--\ref{thm:31oy13thqog} can be translated into a pictorial representation using the idea of motifs (see Figure \ref{fig:o45yuohdaguh3iur}). 
There are 8 distinct motifs, each has two reactions depicted with red arrows. 
The reactants are connected by a green line segment, the {\em reactant polytope} of the reaction network. See \cite{joshi2017small} where this representation for networks with two reactions was used to identify the multistationarity property. 
An upward pointing arrow indicates a reaction of the type $mY \to nY$ for some $n > m$. An arrow pointing towards the north-west indicates a reaction of the type $a_1 X + b_1 Y \to a_2 X + b_2 Y$ with $a_2 < a_1$ and $b_2 > b_1$, and so on. 

\section{Dynamic ACR in networks with an invariant hyperplane} \label{sec:dynamicacr_classification}

We now consider networks with 2  reactions and 2 or fewer species and allow the stoichiometric space to have dimension either 1 or 2. 
With dimension 2 or higher, there is a possibility of oscillatory solutions. 
However, it is easy to see that at least three reactions are required for such solutions (see also \cite{schnakenberg1979simple}). 
It is worth mentioning in passing that for general reaction networks, oscillations are compatible with both static and dynamic ACR. 
An example of a mass action system with three reactions that has {\em oscillations and static ACR in both species} is the classic Lotka-Volterra system ($A \to 2A, A+B \to 2B, B \to 0$), see \cite{joshi2021foundations} for a discussion on this.
The same paper \cite{joshi2021foundations} also features an example with oscillations and dynamic ACR, but it requires 6 reactions and a three dimensional stoichiometric subspace.
When a network with 2 reactions and 2 species has static ACR, every point on the hyperplane $\HH[x,x^*]$ is a steady state. 
A natural generalization of this property is to consider the family of networks for which the hyperplane $\HH[x,x^*]$ is invariant. 
For $\dim(\SS) = 2$ networks, static ACR is ruled out by results from the previous section but dynamic ACR is still possible.

\begin{theorem} \label{thm:oerigheogh}
Let $(\GG,K)$ be as in \eqref{eq:qpeogoe0thg}-\eqref{eq:34oyhhfiwjfp2j}. 
\been
\item There is a unique $x^* > 0$ such that $\HH[x,x^*]$ is invariant  if and only if $b_2=b_1$, $a_2 \ne a_1$ and $(\wt a_1 - a_1)(\wt a_2 - a_2) < 0$. 
\item There is a unique $x^* > 0$ such that $\HH[x,x^*]$ is globally weakly attracting if and only if $\HH[x,x^*]$ is globally weakly stable if and only if $b_2=b_1$, $a_2 \ne a_1$, $(\wt a_1 - a_1)(\wt a_2 - a_2) < 0$ and $(a_2-a_1)(\wt a_1 - a_1) > 0$. 
\enen
\end{theorem}
\begin{proof}
The condition that $0 = \dot x|_{x=x^*}$ is equivalent to: 
\begin{align} \label{eq:qoeghegoh}
0 = k_1(\wt a_1 - a_1) + k_2 (\wt a_2 -a_2) (x^*)^{a_2-a_1} y^{b_2-b_1}. 
\end{align} 
For any $x=x^*$, the equation is an identity if $\wt a_1 = a_1$ and $\wt a_2 = a_2$. So assume that this is not the case. It is clear that there is a unique positive $x^*$ for which \eqref{eq:qoeghegoh} holds if and only if $a_2 \ne a_1$, $b_2 = b_1$ and $(\wt a_1 - a_1)(\wt a_2 - a_2) < 0$. 

For the second part, it suffices to assume that $a_2 \ne a_1$, $b_2 = b_1$ and $(\wt a_1 - a_1)(\wt a_2 - a_2) < 0$ because otherwise $\HH[x,x^*]$ is not invariant and so is not globally weakly attracting. Further we may assume that $a_2>a_1$, possibly after relabeling of reactions. With these assumptions, the original dynamical system $\DD$ in \eqref{eq:34oyhhfiwjfp2j} has the same trajectories as the following dynamical system $\DD'$: 
\begin{align}  \label{eq:wogh3irhglsgn}
\dot x &= k_1(\wt a_1 - a_1)  + k_2(\wt a_2 - a_2) x^{a_2-a_1}  \nonumber \\ 
\dot y &= k_1(\wt b_1 - b_1)  + k_2(\wt b_2 - b_2) x^{a_2-a_1} . 
\end{align}
Note that the $x$ equation is autonomous and does not have a $y$-dependence. 
So a solution to the $x$ equation exists for all time $t \in [0, \infty)$. 
If $\wt a_1 < a_1$ and $\wt a_2 > a_2$ then the $x$-component of the trajectories moves away from $x^*$ monotonically and so $x^*$ is not weakly attracting. 
If $\wt a_1 > a_1$ and $\wt a_2 < a_2$ then the $x$-component of the trajectories moves towards $x^*$. 
This proves the second part. 
\end{proof}
\begin{theorem} \label{thm:qotih3ohoeghow}
Let $(\GG,K)$ be as in \eqref{eq:qpeogoe0thg}-\eqref{eq:34oyhhfiwjfp2j}. Suppose that $\HH[x,x^*]$ is globally weakly attracting for some $x^* > 0$. Denote by $\sigma_i$ the slope of the reaction vector $v_i$, $i \in \{1,2\}$. 
\been
\item $x(t) \not \to \HH[x,x^*]$ for any $(x(0),y(0)) \not \in \HH[x,x^*]$ if $(a_2 - a_1)(\sigma_2-\sigma_1) < 0$. 
\item $x(t) \to \HH[x,x^*]$ for every $(x(0),y(0)) \in \HH[x,x^*] + \spn\{v_2\}$ if $(a_2 - a_1)(\sigma_2-\sigma_1) \ge 0$. 
\item $x(t) \to \HH[x,x^*]$ for every initial value in some cylindrical neighborhood of $\HH[x,x^*]$ if $(a_2 - a_1)(\sigma_2-\sigma_1) > 0$. 
\item $\HH[x,x^*]$ is globally attracting if and only if $\wt b_1 \ge b_1$ and $\wt b_2 \ge b_2$.  
\enen
\end{theorem}
\begin{proof}
From the proof of Theorem \ref{thm:oerigheogh}, the dynamical system $\DD$ is equivalent to $\DD'$ in \eqref{eq:wogh3irhglsgn} with $a_2 > a_1$, $\wt a_1 > a_1$ and $\wt a_2 < a_2$. 
Note that for $\DD'$, $x \to x^*$ and $\dot x \to 0$ asymptotically.
In particular, $\abs{\dot y}$ is bounded for every initial value. 
This means that the solution to the dynamical system $\DD'$ exists for all nonnegative times for every initial value in $\R^2_{\ge  0}$. 
 
Now, the $\dot y$ equation may be written as
\begin{align}
\dot y = k_1(\wt a_1 - a_1) (\sigma_2 - \sigma_1) + \sigma_2 \dot x
\end{align}
If $\sigma_2 = \sigma_1$, then every trajectory converges to some steady state, and this positive steady state is on $\HH[x,x^*]$ when the initial value is compatible with $\HH[x,x^*]$. 
If  $\sigma_2 \ne \sigma_1$, then after some finite time $\dot y$ is bounded away from $0$ and $\sgn(\dot y) = \sgn(\sigma_2 - \sigma_1)$. In particular, for $\sigma_2 < \sigma_1$, $y(t)$ reaches $0$ in finite time. Therefore, the trajectory fails to converge to $\HH[x,x^*]$ for any positive initial value not on $\HH[x,x^*]$. 
On the other hand, if $\sigma_2 > \sigma_1$, then $\dot y > \sigma_2 \dot x$ everywhere in $\R^2_{> 0}$. 
When $\dot y = \sigma_2 \dot x$, we have convergence to $\HH[x,x^*]$ for every $(x(0),y(0)) \in \HH[x,x^*] + \spn\{v_2\}$. So it follows that when $\dot y > \sigma_2 \dot x$, after a finite time, the trajectory enters a cylinder $\CC_\ve$ with $\ve$ such that $\dot y > k_1(\wt a_1 - a_1) (\sigma_2 - \sigma_1)/2$ within $\CC_\ve$. 
Moreover, every trajectory with an initial value in $\CC_\ve$ must converge to the hyperplane $\HH[x,x^*]$. This completes the proof of points {\it 1-3}. 

If $\wt b_1 \ge b_1$ and $\wt b_2 \ge b_2$, then it is immediate from \eqref{eq:wogh3irhglsgn} that $\dot y \ge  0$ everywhere in $\R^2_{\ge 0}$ and so $\HH[x,x^*]$ is globally attracting. 
To show the converse, suppose first that $\wt b_1 < b_1$.  There is a $\delta > 0$ such that for $\{(x(0), y(0)) \in \R^2_{>0}: x(0)^2 + y(0)^2 < \delta^2\}$, $\dot y$ is negative and bounded away from zero. So such a trajectory reaches the $y=0$ boundary and fails to converge to $\HH[x,x^*]$. A similar argument works for $\wt b_2 < b_2$, where we can choose an initial $x(0)$ sufficiently large and an initial $y(0)$ sufficiently small. Thus in either case $\HH[x,x^*]$ is not globally attracting. 
\end{proof}

\begin{theorem}
Let $\GG$ be as in \eqref{eq:qpeogoe0thg}-\eqref{eq:34oyhhfiwjfp2j}.  
Suppose that $X$ is subspace dynamic ACR but not full basin dynamic ACR. Then the following hold:
\been
\item $X$ is a wide basin dynamic ACR species if and only if $(\wt b_i - b_i)(a_i - a_j) > 0$ for $i \ne j \in \{1,2\}$. 
\item $X$ is a narrow basin dynamic ACR species if and only if $(\wt b_i - b_i)(a_i - a_j) < 0$ for $i \ne j \in \{1,2\}$.  
\enen 
\end{theorem}

\begin{proof}
Since $X$ is subspace dynamic ACR, $(a_2 - a_1)(\sigma_2-\sigma_1) \ge 0$ by the previous theorem. We may assume that $a_2 > a_1$, $\sigma_2 \ge \sigma_1$, $\wt a_1 > a_1$ and $\wt a_2 < a_2$. The condition in the first statement is equivalent to $\wt b_1 < b_1$ and $\wt b_2 > b_2$, i.e. $\sigma_2 < 0$. This means that the  set of points incompatible with $\HH[x,x^*]$ is to the left of the hyperplane and  therefore the condition is equivalent to wide basin dynamic ACR. A similar argument in the second case shows equivalence with $\sigma_2 > 0$ and therefore with narrow basin dynamic ACR. 
\end{proof}

\section{Summary of ACR properties for networks with 2 reactions and at most 2 species} \label{sec:summaryofallacr}
The following theorem is a summary of the main results on the different ACR properties in networks with 2 reactions and at most 2 species. 
\begin{theorem} \label{thm:mainesttheorem}
Let $\GG$ be a reaction network with 2 species $\{X, Y\}$ and the following two reactions:
\begin{align} \label{eq:rpyjihrhineoighj3}
a_1 X + b_1 Y \xrightarrow{k_1} \wt a_1 X + \wt b_1 Y, \quad  
a_2 X + b_2 Y \xrightarrow{k_2} \wt a_2 X + \wt b_2 Y, 
\end{align}
where $a_i, b_i, \wt a_i, \wt b_i \in \R_{\ge 0}$ and $(\wt a_i, \wt b_i) \ne (a_i, b_i)$ for  $i \in \{1,2\}$. 
The labels $k_1$ and $k_2$ are mass action reaction rate constants, when considering a mass action system $(\GG,K)$. 

\been

\item A  necessary and sufficient condition for the existence of a unique invariant hyperplane $\HH$ parallel to a coordinate axis is that the reactant polytope be a line segment parallel to the other coordinate axis. In particular, 
$(a_1,b_1) \ne (a_2,b_2)$ and $(a_2-a_1)(b_2 - b_1)=0$. 

\suspend{enumerate}

\noindent Suppose that $\GG$ has a unique invariant hyperplane $\HH$ parallel to a coordinate axis. Then the following hold: 

\resume{enumerate}

\item $\GG$ has the capacity for $\PP$-ACR if and only if $\GG$ is $\PP$-ACR, where $\PP \in \{$static, strong static, dynamic, weak dynamic$\}$, i.e. the $\PP$-ACR properties are independent of choice of rate constants.

\item We use  the convention $0\cdot (1/0)=0$ in the following expressions. 
\been
\item $\GG$ is static ACR if and only if $\GG$ is strong static ACR if and only if  $\dim(\SS) = 1$ and reaction vectors are negative scalar multiples of each other:
$\wt a_1-a_1 = - \mu \left(\wt a_2-a_2\right) \mbox{ for some } \mu > 0.$
\item $\HH$ is weakly attracting if and only if both reactions point inwards:
\begin{align*}
(\wt a_i - a_i) (a_j - a_i) + (\wt b_i - b_i) (b_j - b_i) > 0, \quad  i \ne j \in \{1,2\}. 
\end{align*}
\item $\GG$ is non-null dynamic ACR if  and only if $\GG$ is subspace dynamic ACR if  and only if 
\begin{align*}
\sum_{i \ne j \in \{1,2\}} \left\{(a_j-a_i)  \left(\frac{\wt b_i - b_i}{\wt a_i - a_i}\right) + (b_j-b_i) \left(\frac{\wt a_i - a_i}{\wt b_i - b_i}\right) \right\} \ge 0. 
\end{align*}

\suspend{enumerate}

\noindent Suppose that neither $\dot x$ nor $\dot y$ is identically zero, or equivalently $(\wt a_1, \wt a_2) \ne (a_1, a_2)$ and $(\wt b_1, \wt b_2) \ne (b_1, b_2)$. 

\resume{enumerate}

\item $\GG$ is cylinder dynamic ACR if  and only if 
\begin{align*}
\sum_{i \ne j \in \{1,2\}} \left\{(a_j-a_i)  \left(\frac{\wt b_i - b_i}{\wt a_i - a_i}\right) + (b_j-b_i) \left(\frac{\wt a_i - a_i}{\wt b_i - b_i}\right) \right\} > 0. 
\end{align*}

\item $\GG$ is full basin dynamic ACR if and only if 
\begin{align*}
(a_j-a_i) \left(\frac{\wt b_i - b_i}{\wt a_i - a_i}\right) + (b_j-b_i) \left(\frac{\wt a_i - a_i}{\wt b_i - b_i}\right) \ge 0, \quad i \ne j \in \{1,2\}. 
\end{align*}
\enen
\item For all the above networks with $\PP$-ACR: 
\beit
\item Either $X$ or $Y$, but not both, is an $\PP$-ACR species. 
\item $X$ is an $\PP$-ACR species if $a_2 \ne a_1$. The variable $x$ has the $\PP$-ACR value $\ds \left( - \frac{k_2 (\wt a_2 - a_2)}{k_1 (\wt a_1 - a_1)} \right)^{\frac{1}{a_1-a_2}}$. 
\item $Y$ is an $\PP$-ACR species if $b_2 \ne b_1$. The variable $y$ has the $\PP$-ACR value $\ds \left( - \frac{k_2 (\wt b_2 - b_2)}{k_1 (\wt b_1 - b_1)} \right)^{\frac{1}{b_1-b_2}}$. 
\enit

\enen

\end{theorem}

A pictorial representation of these results is shown in Figure  \ref{fig:p318tyhelgheg}.

\section{Discussion and future work}

In this paper, we have established that for small networks, the Euclidean embedding (or geometric structure) of a reaction network can yield deep insights into the dynamics of the mass action ODE system. In particular, when there are only two reactions and at most two species, the reactant polytope (the line segment joining the two reactant complexes) is required to be horizontal or vertical for any type of ACR property. 
Some of the two reaction networks that appear in this paper do show up in applications. For instance, the archetypal wide basin ACR network in \ref{fig:qperuhg305y84oigg}(a) can be thought of as a simple model of infectious disease dynamics (SIS model). The reaction $A+B \to 2B$ represents an infective individual $B$ infecting a susceptible individual $A$, while the reaction $B \to A$ represents recovery from infection. 
One may also interpret this model as a protein with two alternate conformations $A$ and $B$ with the  spontaneous transition $B \to A$, while the transition $A \to B$ is catalyzed (or promoted) by $B$. 
However, most biochemically realistic networks are significantly more complex with several reactions and species.  
Even though ACR can be found in higher dimensional systems, we do not expect that there will be such simple characterization of their ACR properties. 
However, our work suggests that it may be fruitful to study the link between the geometry of reactant polytopes and dynamics of ACR systems further. 
Moreover, we do expect that small motifs that are embedded within large and complicated networks may affect the overall dynamics. 
Our future goal is to understand such effects. 

The small motifs studied in this paper also serve as test cases for various dynamical behaviors. 
For instance, a surprising possibility revealed from the study of small motifs was that of weak dynamic ACR, where every trajectory monotonically approaches a hyperplane while simultaneously failing to converge (see Figure \ref{fig:o;eithjeoihy3j2potj}(b)). 

In future work, we study dynamic ACR in more complex, biochemically realistic systems such as bacterial two-component signaling systems, a class that encompasses several thousands of systems \cite{alon2019introduction,barkai1997robustness,alon1999robustness,batchelor2003robustness,shinar2007input,shinar2009robustness}. 
Moreover, we study {\em consequences} of dynamic ACR. We plan to show that dynamic ACR, i.e. the property of robustness against variations in initial conditions, surprisingly leads to other much more robust dynamical properties with stronger implications for biochemical systems with dynamic ACR.

\newpage
\subsection*{Acknowledgments}
G.C. was supported by NSF grant DMS-2051568 and by a Simons Foundation fellowship. 
We thank the referees for careful reading and helpful comments.
\bibliographystyle{plain}
\bibliography{acr}

\end{document}